\documentclass[12pt]{amsart}
\usepackage[dvips]{graphicx}

\begin{document}

\def\Z{{\mathcal Z}}
\def\V{{\mathcal V}}
\def\F{{\mathcal F}}
\def\I{{\mathcal I}}

\def\affine{\mathbf{A}}
\def\variety{\Z^{m,n}_{r,k}}
\def\lvariety{\Z^{m-1,n-1}_{r-1,k}}
\def\sqvariety{\Z^{m,m}_{m,k}}
\def\gensqvariety{\Z^{m,m}_{r,k}}
\def\ideal{\I^{m,n}_{r,k}}
\def\sqideal{\I^{m,m}_{m,k}}
\def\lideal{\I^{m-1,n-1}_{r-1,k}}

\def\p{\mathbf{p}}
\def\u{\mathbf{u}}
\def\w{\mathbf{w}}
\def\x{\mathbf{x}}
\def\y{\mathbf{y}}
\def\vec{\mathbf{v}}
\def\A{{\mathbf A}}
\def\P{ {\mathbf P} }
\def\Aff{\A^{n^2} }
\def\gen{{\mathbf X}}

\def\pf{\medskip\noindent{\bf Proof.} \ }
\def\qed{{\hfill{$\Box $}}\vskip 6pt}
\def\C{{\mathcal C}}
\def\CP{{\mathcal C_{\bf P}}}

\def\trunc{F[t]/(t^k)}
\def\vars{x^{(l)}_{i,j}}
\def\zvars{z^{(l)}_{i,j}}

\def\Z{{\mathcal Z}}
\def\V{{\mathcal V}}
\def\F{{\mathcal F}}
\def\Gb{{\mathcal G}}
\def\G{{{\mathbf G}(n,n^2)}}
\def\text{}
\def\im{\mathrm{image}}
\def\hY{\widehat{Y}(i,j)}
\def\hYkl{\widehat{Y}(k,l)}
\def\tY{\widetilde{Y}(i,j)}
\def\tYkl{\widetilde{Y}(k,l)}
\def\lm{\textsc{lm}}

\newtheorem{Theorem}{Theorem}[section]
\newtheorem{Lemma}[Theorem]{Lemma}
\newtheorem{Proposition}[Theorem]{Proposition}
\newtheorem{Corollary}[Theorem]{Corollary}

\theoremstyle{remark}
\newtheorem{Remark}[Theorem]{Remark}
\newtheorem{Notation}[Theorem]{Notation}
\newtheorem{Definition}[Theorem]{Definition}
\newtheorem{Note}[Theorem]{Note}
\newtheorem{Example}[Theorem]{Example}

\author {Toma\v{z} Ko\v{s}ir}
\author{B.A. Sethuraman}
\address{Dept. of Mathematics\\University of Ljubljana\\Jadranska
19\\ 1000 Ljubljana\\Slovenia}
\address{Dept. of Mathematics\\California State University
Northridge\\Northridge CA 91330\\U.S.A.}
\title{Determinantal Varieties Over Truncated Polynomial Rings}
 \email{tomaz.kosir@fmf.uni-lj.si}
\email{al.sethuraman@csun.edu}

\begin{center}
\today
\end{center}

\begin{abstract}
We study higher order determinantal varieties obtained by
considering generic $m\times n$ ($m \le n$) matrices over rings
of the form $\trunc$,  and for some fixed $r$, setting the
coefficients of powers of $t$ of all $r \times r$ minors to zero.
These varieties can be interpreted as generalized tangent bundles
over the classical determinantal varieties; a special case of
these varieties first appeared in a problem in commuting
matrices. We show that when $r = m$, the varieties are
irreducible, but when $r < m$, these varieties have at least
$\lfloor {k/2}\rfloor + 1$ components. In fact, when $r=2$ (for
any $k$), or when $k=2$ (for any $r$), there are exactly $\lfloor
{k/2}\rfloor + 1$ components. We give formulas for the dimensions
of these components in terms of $k$, $m$, and $n$.  In the case
of square matrices with $r=m$, we show that the ideals of our
varieties are prime and that the coordinate rings are complete
intersection rings, and we compute the degree of our varieties
via the combinatorics of a suitable simplicial complex.

\end{abstract}
\maketitle

\section{Introduction} Let $F$ be an algebraically closed field
and $\affine_F^{k}$ the affine space of dimension $k$ over $F$.
 The varieties $\Z^{m,n}_r \subset
\affine_F^{mn}$ consisting of $m\times n$ matrices ($m \le n$)
with entries in $F$ and of rank at most $r-1$ are of course a
natural and  very well understood class of objects; their various
geometric and algebraic properties and their connections to
representation theory and combinatorics have been extensively
studied (see \cite{BV} for instance).  By contrast, very little
is known about the following class of objects $\Z^{m,n}_{r,k}$
that are very closely related to the classical varieties
$\Z^{m,n}_r$: Consider the truncated polynomial ring
$R=F[t]/(t^k)$ ($k = 1, 2,3,\dots$), and let $X(t)$ be the
generic $m\times n$ matrix over this ring; thus, each entry of
$X$ is of the form  $ x_{i,j}(t) = x_{i,j}^{(0)} + x_{i,j}^{(1)}
t + \cdots + x_{i,j}^{(k-1)} t^{k-1}$.
  Each $r \times r$
minor of this matrix is an element of $R=F[t]/(t^k)$. Let
$\I^{m,n}_{r,k}$  be the ideal of $F[\{x_{i,j}^{(l)},\ 1 \le i
\le m,\ 1 \le j \le n,\ 0 \le l < k\}]$ generated by the
coefficients of $t$ in each $r \times r$ minor of the generic
matrix $X(t)$, and define $\Z^{m,n}_{r,k} \subseteq
\affine_F^{nmk}$ to be the zero set of $\I^{m,n}_{r,k}$.  These
varieties $\Z^{m,n}_{r,k}$ are therefore natural generalizations
of the classical varieties $\Z^{m,n}_{r}$, and when $k=1$, of
course, we simply recover the original $\Z^{m,n}_{r}$.

Our interest in these varieties arises from previous work  on
commuting triples of matrices.  In the paper \cite{NS}, the
second author and Neubauer determined the variety of commuting
pairs in the centralizers of $2$-regular matrices (a matrix is
said to be $r$-regular if each eigenspace is at most $r$
dimensional).  They observed there that when $C$ is a $2$-regular
$n\times n$ matrix, the variety of commuting pairs in the
centralizer of $C$ is the product of $\affine_F^p$ (for suitable
$p$) and the subvariety of $2 \times 3$ matrices over
$F[t]/(t^k)$ where the coefficients of $t$ of all $2\times 2$
minors vanish. This second factor is of course just the variety
$\Z^{2,3}_{2,k}$ introduced above.  It was then natural to
recognize $\Z^{2,3}_{2,k}$ as belonging to the larger class of
varieties $\Z^{m,n}_{r,k}$, and to commence a program to study
this larger family.

It is worth giving a geometric interpretation to these varieties
$\Z^{m,n}_{r,k}$ as suitable ``bundles'' over the classical
objects $\Z^{m,n}_{r}$.  Suppose $f(x,y,\dots)$ is an irreducible
polynomial, and suppose we were to replace the variables $x$ by
$\sum_0^{k-1} x^{(i)}t^i$, $y$ by $\sum_0^{k-1} y^{(i)}t^i$, and
so on, and suppose we were to expand the polynomial in powers of
$t$ and write it as $f_0 + f_1 t + f_2 t^2 + \cdots$.  (Note that
$f_0$ will just be $f$.) Setting each of $f_0$, $f_1$, $f_2$,
$\dots$, $f_{k-1}$ to zero is to ask for a point $(x^{(0)},
y^{(0)}, \dots)$ on the hypersurface $f=0$, and then to ask for
those degree $k-1$ curves $x = x^{(0)}+\sum_1^{k-1} x^{(i)}t^i$,
$y = y^{(0)}+\sum_1^{k-1} y^{(i)}t^i$, and so on, parameterized
by $t$, that vanish   on the variety of $f$ at that point up to
order $k$. Applying these considerations to the equations
defining our varieties $\Z^{m,n}_{r,k}$ (and recalling that the
obvious defining equations for $\Z^{m,n}_{r}$ also generate the
ideal of $\Z^{m,n}_{r}$), it is clear then that $\Z^{m,n}_{r,k}$
consists of the classical varieties $\Z^{m,n}_{r}$, and at each
point of $\Z^{m,n}_{r}$, those parameterized degree $k-1$ curves
vanishing on $\Z^{m,n}_{r}$ at that point to order $k$. In
particular, when $k=2$, $\Z^{m,n}_{r,2}$ may be considered as the
''tangent bundle'' to the classical determinantal variety
$\Z^{m,n}_{r}$.  Of course, this is not really a bundle in the
usual sense, since different fibers will have different
dimensions, but we will use the words tangent bundle freely in
the paper.  (More accurately, $\Z^{m,n}_{r,2}$ is the variety
associated to the symmetric algebra on the module of derivations
of the classical determinantal variety $\Z^{m,n}_{r}$.)

In a different language, the varieties $\variety$ appear  as the
$(k-1)$-th jet scheme of the classical determinantal varieties
$\Z_{r}^{m,n}$ \cite[\S 1]{MM}. From another point of view, the
varieties $\variety$ are the restriction (or direct image) from
$R$ to $F$, in the sense of Weil (\cite[I, \S 1, 6.6]{DG}) of the
scalar extension of the classical varieties $\Z^{m,n}_{r} \times_F
R$.

Since in general   the fibers over the base $\Z^{m,n}_{r}$ will
not all be of the same dimension, it is not a priori clear
whether the assemblage of the base space and its fibers should be
reducible or irreducible.  We show here that if we set all
maximal minors to zero (that is, if we take $r = m$), then
$\Z^{m,n}_{m,k}$ is indeed irreducible, but if we consider
\textit{submaximal} minors (that is, if we take $r \le m-1$),
then $\Z^{m,n}_{r,k}$ breaks up into several components, not all
of the same dimension.  In fact, we have a complete picture of
what these various components look like in the case when $r = 2$,
i.e., in the case of $2\times 2$ minors (for any $k$); we also
have a complete picture of the components   when $k=2$ (for any
$r$). In both these situations, there are exactly $\lfloor
{k/2}\rfloor + 1$ components (so in particular the tangent bundle
has two components), and we have formulas for their dimensions in
terms of $k$, $m$, and $n$.  In general, we show that when $r <
m$, the minimum number of components must be $\lfloor
{k/2}\rfloor + 1$.  Also, from the fact that these components
intersect nontrivially, we are able to show that in the
submaximal minors case, our varieties are not normal, and that
since the components (with one exception) are of different
dimensions, that our varieties (with one possible exception) are
not Cohen-Macaulay.

In the special case where $m=n$, and where we consider maximal
minors, we are able to say considerably more.  We show in that
case that the defining equations form a Groebner basis (with
respect to a suitable ordering) for the ideal generated by these
equations, and we are then able to show that the ideal is
actually prime. It follows easily that the coordinate ring is a
complete intersection ring.  Moreover, the quotient of the
polynomial ring in $m^2k$ variables by  the ideal generated by
the leading terms of these equations turns out to be the
Stanley-Reisner ring of a particularly nice simplicial complex,
and from the combinatorics of that complex, we are able to
determine the Hilbert polynomial of our original ideal
$\I^{m,m}_{m,k}$.

We introduce here some alternative notation for the entries of
the generic matrix that will be of much use: We will  denote
  the $i$-th row of the matrix $X(t)$ over
$\trunc$ by $\u_{i}(t)$: this is an element of $(F[t]/t^k)^n$. We
will write $\u_{i}(t) = \sum_{l=0}^{k-1} \u_i^{(l)} t^l$, so the
various $\u_i^{(l)}$ are row vectors from $F^n$.  We will
sometimes refer to $\u_i^{(l)}$ by itself as the ``row
$\u_i^{(l)}$.''  We will also refer to a vector of the form
$\u_i^{(l)}$ as being ``of degree $l$.'' In a similar vein, we
will talk of a variable of the form $\vars$ as being of ``of
degree $l$.''  In particular, when we talk of a ``degree zero''
minor, we will mean a minor of the matrix $X(0) =
((x^{(0)}_{i,j}))$.

The methods we use in the paper are totally elementary. We also
note that the paper \cite{SoSt} discusses a related set of
objects: the quantum Grassmannians, whose coordinate rings are
the subalgebras of $F[\vars]$ generated by the coefficients of
$t$ of the various $m\times m$ minors of an $m\times n$ matrix.

\section{The Fundamental Reduction Process} \label{reductionsetup}

We describe here a reduction process that exhibits our varieties
$\Z^{m,n}_{r,k}$ to be a union of two subvarieties, one isomorphic
to $\Z^{m,n}_{r,k-r} \times \affine^{mn(r-1)}$ (or to
$\affine^{mn(k-1)}$ when $k \le r$), and another whose components
are in one-to-one correspondence with the components of
$\lvariety$, and which, in fact, is birational to $\lvariety
\times \affine^{(m+n-1)k}$.  This reduction will be a key tool in
understanding the components of $\variety$.

\begin{Lemma} \label{allbasevarszero} The subvariety of $\Z^{m,n}_{r,k}$ where all
$x^{(0)}_{i,j}$ are zero is isomorphic to $\Z^{m,n}_{r,k-r} \times
\affine^{mn(r-1)} $ when $k > r$, and isomorphic to
$\affine^{mn(k-1)}$ when $k \le r$.
\end{Lemma}
\begin{proof}  This can be seen easily by writing the equations
defining $\Z^{m,n}_{r,k}$ in terms of the rows $\u_i^{(l)}$.  Our
determinantal equations read
\begin{equation}
\u_{i_1}\wedge \u_{i_2}\wedge\dots\wedge\u_{i_r} = \mathbf{0}
\end{equation}
for all $1 \le i_1 < i_2 < \dots < i_r \le m$.  This is an
equation in $\bigwedge^r (\trunc)^n$, which expands to  a set of
$k$ equations in $\bigwedge^r F^n$, one for the coefficient of
each power of $t$. The equation for the coefficient of $t^l$ reads
\begin{equation}
\sum_{d_1+d_2+\cdots+d_r = l}\u_{i_1}^{(d_1)}\wedge
\u_{i_2}^{(d_2)}\wedge\dots\wedge\u_{i_r}^{(d_r)} =
\mathbf{0},\quad l = 0, \dots, k-1.
\end{equation}
It is clear that if all $\u_i^{(0)}$ are zero, then all terms in
the coefficients of $t^l$, for $l = 0, 1, \dots, r-1$ become
zero, since any $r$-fold product of degree at most $r-1$ must
contain at least one term of degree $0$. If $k \le r$, this just
means that there are no equations governing the remaining
variables $x^{(l)}_{i,j}$ for $l\ge 1$, so the subvariety is
isomorphic to $\affine^{mn(k-1)}$. When $k >r$, the equation for
the coefficient of $t^l$ for $l = r, r+1, \dots, k-1$ now reads
(after all terms involving any $\u_i^{(0)}$ have been removed)
\begin{equation}
\sum_{\stackrel{d_1+d_2+\cdots+d_r = l}{d_1\ge 1, \dots, d_r\ge 1}
}\u_{i_1}^{(d_1)}\wedge
\u_{i_2}^{(d_2)}\wedge\dots\wedge\u_{i_r}^{(d_r)} =
\mathbf{0},\quad l = r, \dots, k-1.
\end{equation}
Observe that none of the rows  $\u_{i}^{(k-1)}$,
$\u_{i}^{(k-2)}$, $\dots$, $\u_{i}^{(k-(r-1))}$ show up in these
equations.  For, every summand is an $r$-fold wedge product of
degree $l$, and if, for instance, $\u_{i}^{(k-(r-1))}$ were to
appear in a summand, then the minimum degree of that summand
would be $k-(r-1) + (r-1) = k > k-1$. Thus, there are no
equations governing the variables $x^{(l)}_{i,j}$ for $k-(r-1)
\le l \le k-1$, which accounts for the factor
$\affine^{mn(r-1)}$. Setting $e_i = d_i-1$, these equations can
be rewritten as
\begin{equation}
\sum_{\stackrel{e_1+e_2+\cdots+e_r = l-r }{e_i \ge 0
}}\u_{i_1}^{(e_1+1)}\wedge
\u_{i_2}^{(e_2+1)}\wedge\dots\wedge\u_{i_r}^{(e_r+1)} =
\mathbf{0}, \quad l = r, \dots, k-1,
\end{equation}
or what is the same thing,
\begin{equation}
\sum_{\stackrel{e_1+e_2+\cdots+e_r = l' }{e_i \ge 0
}}\u_{i_1}^{(e_1+1)}\wedge
\u_{i_2}^{(e_2+1)}\wedge\dots\wedge\u_{i_r}^{(e_r+1)} =
\mathbf{0}, \quad l' = 0, \dots, k-r-1.
\end{equation}
 But these are precisely the equations
that one would obtain if one were to consider the generic matrix
$m\times n$ matrix with rows $\u_{i}^{(1)} + \u_{i}^{(2)} t +
\cdots + \u_{i}^{(k-r)} t^{k-r-1}$ ($1 \le i \le m$) and set
determinants of $r \times r$ minors to zero modulo $t^{k-r}$.
This proves the lemma.
\end{proof}

Our next theorem will be crucial to understanding the closure of
the open set where at least one $\vars$ is nonzero. It is merely
an extension to the case $k > 1$ of a result that is well known
in the classical case (see \cite[Prop. 2.4]{BV}, for instance).

We first need the following:
\begin{Lemma} \label{anygensOK}
Let $R$ be a ring, and let $I$ be an ideal of $R[t]/(t^k)$. Let
$I$ be generated by $a_1(t), \dots, a_u(t)$, with each $a_i =
\sum_{0}^{k-1} a_i^{(l)} t^l$.  Let $J$ be the ideal of $R$
generated by the various $a_i^{(l)}$.  Then $J=K$, where $K$ is
the set of all $r\in R$ such that $r$ is the coefficient of
$t^l$, for some $l$, in some element of $I$.  In particular, if
$I$ is also generated by $b_1(t),\dots, b_v(t)$ with $b_i =
\sum_{0}^{k-1} b_i^{(l)} t^l$, then the ideal of $R$ generated by
the various $b_i^{(l)}$ also equals $J$.
\end{Lemma}
\begin{proof} This is elementary.
\end{proof}

Write $S = F[\vars]$, and observe that in the ring
$S[({x^{(0)}_{m,n}})^{-1}][t]/(t^k)$, the element $x_{m,n}(t)$ is
invertible. If we were to perform row reduction in
$S[({x^{(0)}_{m,n}})^{-1}][t]/(t^k)$
on the matrix $X$ to bring all the elements in the last column
above   $x_{m,n}$ to zero, we would subtract from row $\u_i$ the
row $\u_m$ multiplied by $x_{m,n}^{-1} x_{i,n}$. Thus, we would
replace $X$ by a matrix
$$ Y = ((y_{i,j})),$$
where
\begin{equation} \label{ymat}
y_{i,j} = \left\{ \begin{array}{ll} x_{i,j} -
x_{m,j}x_{i,n}x_{m,n}^{-1} & \mbox{ for $1 \le i \le m-1,\ 1 \le j \le n-1$} \\
0 & \mbox{for $j = n$ and $1 \le i \le m-1$}\\
x_{i,j} & \mbox{for $i = m$ and $1 \le j \le n$}
\end{array}
\right.
\end{equation}

Since the inverse of $x_{m,n}(t)$ can be written as a polynomial
in the various entries $x^{(l)}_{m,n}$ (for $1 \le l \le k-1$)
and various negative powers of $x^{(0)}_{m,n}$, we find that each
$y^{(l)}_{i,j}$, for $1 \le i \le m-1,\ 1\le j \le n-1$ can be
written in terms of the $\vars$ in the form
\begin{equation} \label{yform}
y^{(l)}_{i,j} = \vars - q^{(l)}_{i,j}(x^{(p)}_{m,j},
x^{(r)}_{i,n}, x^{(s)}_{m,n}, (x^{(0)}_{m,n})^{-1}),
\end{equation}
for a suitable polynomial expression $q^{(l)}_{i,j}$ in the
indicated variables ($0 \le p, r < k$, $1 \le s < k$).

With the observations above about row reduction as our
motivation, and continuing with   the same notation,   we have
the following:

\begin{Theorem} (see \cite[Prop. 2.4]{BV}) \label{induct}  Assume $r \ge 2$.  Let $\zvars$, $1 \le i
\le m-1,\ 1 \le j \le n-1$, $0 \le l < k$ be a new set of
variables, and write $T$ for the ring $F[\zvars]$, and $T'$ for
the ring $F[\zvars,
x^{(l)}_{1,n},\dots,x^{(l)}_{m,n},x^{(l)}_{m,1},\dots,x^{(l)}_{m,n-1}]$
($0 \le l < k $). Also, write $Z$ for the $m-1 \times n-1$ matrix
$((z_{i,j}(t)))$ over  $T[t]/(t^k)$,  where $z_{i,j}(t) =
\sum_{l=0}^{k-1} \zvars t^l$.   We have an isomorphism
$$
S[(x^{(0)}_{m,n})^{-1}] \cong T'[(x^{(0)}_{m,n})^{-1}],
$$ given by
\begin{eqnarray*}
f\colon \quad
x^{(l)}_{i,n} &\rightarrow & x^{(l)}_{i,n}\ 1 \le i \le m\\
x^{(l)}_{m,j} &\rightarrow & x^{(l)}_{m,j}\ 1 \le j \le n-1\\
\vars &\rightarrow & \zvars + q^{(l)}_{i,j}(x^{(p)}_{m,j},
x^{(r)}_{i,n}, x^{(s)}_{m,n},
(x^{(0)}_{m,n})^{-1}), \ 0 \le p,r < k, \ 1\le s < k,\\
{} && \mbox{for $1 \le i \le m-1,\ 1\le j \le n-1$, $0 \le l <
k$.}
\end{eqnarray*}

Under this isomorphism, the localization of $\ideal$ at
$(x^{(0)}_{m,n})^{-1}$ corresponds to the localization of the
ideal $\lideal T'$  at $(x^{(0)}_{m,n})^{-1}$, where $\lideal$ is
the ideal of $T$ determined by the coefficients of $t$ of the
various  $(r-1) \times (r-1)$   minors of the matrix $Z$.
Moreover, this gives a one-to-one correspondence between the
prime ideals $P$ of $S$ that are minimal over $\ideal$ and do not
contain $x^{(0)}_{m,n}$ and the prime ideals $Q$ of $T$ that are
minimal over $\lideal.$ In this correspondence, if $P$
corresponds to $Q$, then the codimension of $P$ in $S$ equals the
codimension of $Q$ in $T$.
\end{Theorem}

\begin{proof}

The fact that $f$ is an isomorphism is clear, since the map
\begin{eqnarray*}
\tilde{f}\colon \quad
x^{(l)}_{i,n} &\rightarrow & x^{(l)}_{i,n}\ 1 \le i \le m\\
x^{(l)}_{m,j} &\rightarrow & x^{(l)}_{m,j}\ 1 \le j \le n-1\\
 \zvars &\rightarrow & \vars -
q^{(l)}_{i,j}(x^{(p)}_{m,j}, x^{(r)}_{i,n}, x^{(s)}_{m,n},
(x^{(0)}_{m,n})^{-1}), \ 0 \le p,r < k, \ 1\le
s < k,\\
{} && \mbox{for $1 \le i \le m-1,\ 1\le j \le n-1$, $0 \le l < k$}
\end{eqnarray*}
provides the necessary inverse.

As for the second assertion, write $I$ for the   localization of
$\ideal$ at $(x^{(0)}_{m,n})^{-1}$ and $J$ for the localization
of the $\lideal T'$  at $(x^{(0)}_{m,n})^{-1}$.  We wish to show
that $f(I) = J$. Subtracting a multiple of the $m$-th row from
the $i$-th row of a matrix preserves the ideal of
$S[(x^{(0)}_{m,n})^{-1}][t]/(t^k)$  generated by $r\times r$
minors, so by Lemma \ref{anygensOK}, the coefficients of $t$ of
the various $r\times r$ minors of the matrix $Y$ can be taken as
the generators of $I$. Write $\tilde{Y} $ for the upper-left
$m-1\times n-1$ block of $Y$. Recall that $r\ge 2$. The $r\times
r$ minors of $Y$ fall into two classes. The first class consists
of minors that involve the last column of $Y$, so by Laplace
expansion, these minors are either zero or of the form
$x_{m,n}(t)$ times an $(r-1)\times(r-1)$ minor of $\tilde{Y}$.
Since $x_{m,n}(t)$ is invertible,  we find that up to
multiplication by a unit, this class of generators is precisely
the set of all $(r-1)\times(r-1)$ minors of $\tilde{Y}$.   The
second class of generators of $I$ consists of minors that do not
involve the last column of $Y$.  By Laplace expansion, these
minors can be written as $S[({x^{(0)}_{m,n}})^{-1}]$ linear
combinations of suitable $(r-1)\times(r-1)$ minors of
$\tilde{Y}$.  It follows that $I$ is generated precisely by  the
set of all $(r-1)\times(r-1)$ minors of $\tilde{Y}$. On the other
hand, $J$ is generated by all $(r-1)\times (r-1)$ minors of the
matrix $Z$. Thus, under the map $f$, these generators of $I$ map
precisely to generators of $J$, so $f(I)= J$.

As for the last assertion, we have a one-to-one correspondence
between the minimal primes of $\ideal$ that do not contain
$x^{(0)}_{m,n}$ and the minimal primes of the localized ideal $I$
in $S[(x^{(0)}_{m,n})^{-1}]$.  By the isomorphism described
above, these are in one-to-one correspondence with the minimal
primes of the ideal $J$ of $T'[(x^{(0)}_{m,n})^{-1}]$. These, in
turn, are in one-to-one correspondence with the minimal primes of
the ideal $\lideal T'$ of $T'$ that do not contain
$x^{(0)}_{m,n}$.  But since $T'$ is just an extension of $T$
obtained by adding the indeterminates $x^{(l)}_{1,n}$, $ \dots,$
$x^{(l)}_{m,n},$ $x^{(l)}_{m,1}$, $\dots$, $x^{(l)}_{m,n-1}$, the
minimal primes of $\lideal T'$ are in one-to-one correspondence
with the minimal primes of $\lideal$ in $T$;  specifically,    the
minimal prime $Q$ of $\lideal$   corresponds to
$Q[x^{(l)}_{1,n},\dots,x^{(l)}_{m,n},x^{(l)}_{m,1},\dots,x^{(l)}_{m,n-1}]$.
  Moreover, tracing through this correspondence, since
localization and adding indeterminates does not change the
codimension of a prime ideal that avoids the localization set, we
find that the correspondence preserves the codimension of the
respective prime ideals in their respective rings.   This gives
us the assertion.
\end{proof}

\begin{Remark}\label{birational} The theorem above shows that
there is a birational isomorphism between $\variety$ and
  $\lvariety \times \affine^{k(m+n-1)}$,
with the domain of definition of this isomorphism being the open
set of $\variety$ where $x^{(0)}_{m,n} \neq 0$,   and the image
being  the open set of $\lvariety \times \affine^{k(m+n-1)}$
  where the ``free'' variable $x^{(0)}_{m,n}
\neq 0$.
\end{Remark}

\begin{Remark}\label{anyvariable}
Notice that there is nothing special in these considerations
about the variable $x^{(0)}_{m,n}$. Essentially the same result
holds for localization at any other variable $x^{(0)}_{i,j}$. The
matrix $\tilde Y$ in that case would arise from the removal of the
$j$-th column and $i$-th row of the matrix $X$.
\end{Remark}

The ideas in the proof of the theorem above also lead to  the
following:

\begin{Proposition}\label{doesntmatterwhichxij} Let $S = F[\vars]$ as before.  Let $P$ be a minimal
prime ideal of $\ideal$.  Then for any two pairs of indices
$(i,j)$ and $(i',j')$, $P$ contains $x^{(0)}_{i,j}$ iff it
contains $x^{(0)}_{i',j'}$.
\end{Proposition}
\begin{proof} It is sufficient to prove this for the case where
$(i',j') = (m,n)$.  Let $P$ be a prime ideal minimal over $\ideal$
that does not contain $x^{(0)}_{m,n}$.  Assume to the contrary
that it contains $x^{(0)}_{i,j}$.  Then the localization $\tilde
P$ of $P$ at $x^{(0)}_{m,n}$  will also contain $x^{(0)}_{i,j}$,
and will be minimal over the localization of $I$ of $\ideal$.
Hence, the ideal $f(\tilde P)$, where $f$ is as in the theorem
above, will contain $f(x^{(0)}_{i,j}) = z^{(0)}_{i,j} +
q^{(0)}_{i,j}$, and will be minimal over $J = f(I)$. But as is
readily seen, $q^{(0)}_{i,j}$ is just
$x^{(0)}_{i,n}x^{(0)}_{m,j}(x^{(0)}_{m,n})^{-1}$.  Since
$x^{(0)}_{m,n}$ is a unit, it follows that $f(P)$ will contain
$x^{(0)}_{m,n}z^{(0)}_{i,j} + x^{(0)}_{i,n}x^{(0)}_{m,j}$.  Under
the localization map from $T'$ to $T'[(x^{(0)}_{m,n})^{-1}]$,
$f(P)$ will correspond to a prime ideal $Q'$ of $T'$,  that is
minimal over the ideal $\lideal  T'$. Moreover, $Q'$ will contain
$x^{(0)}_{m,n}z^{(0)}_{i,j} + x^{(0)}_{i,n}x^{(0)}_{m,j}$. As we
saw in the proof of the last assertion of Theorem \ref{induct}
above, $Q'$ must be of the form
$Q[x^{(l)}_{1,n},\dots,x^{(l)}_{m,n},x^{(l)}_{m,1},\dots,x^{(l)}_{m,n-1}]$
for some minimal prime $Q$ of the ideal $\lideal$ of $T$. But
this is impossible, since the $Q$ coefficient of
$x^{(0)}_{i,n}x^{(0)}_{m,j}$ in the element
$x^{(0)}_{m,n}z^{(0)}_{i,j} + x^{(0)}_{i,n}x^{(0)}_{m,j}$,
namely, $1$, is not in $Q$.  Hence $P$ cannot contain
$x^{(0)}_{i,j}$.

Conversely, if $P$ is a minimal prime ideal of $\ideal$ that does
not contain $x^{(0)}_{i,j}$ but contains $x^{(0)}_{m,n}$, then
the same argument as above, applied to the corresponding
isomorphism obtained on localizing at $x^{(0)}_{i,j}$ (see Remark
(\ref{anyvariable}) above), gives us a contradiction.  This proves
the corollary.
\end{proof}

We are now ready to decompose our variety $\variety$ as described
at the beginning of this section. Let $Z_0$ represent the union
of the zero sets of all those minimal prime ideals of $\ideal$ in
$S = F[\vars]$ that do not contain  some $x^{(0)}_{i,j}$ (and
hence, by Proposition \ref {doesntmatterwhichxij} above, do not
contain any $x^{(0)}_{i,j}$ for $1 \le i \le m$, $1 \le j \le
n$).  $Z_0$ is not empty: There are clearly points in our variety
where $x^{(0)}_{i,j}$ is not zero, hence there exist minimal
primes of $\ideal$ that do not contain  $x^{(0)}_{i,j}$ (for any
$(i,j)$).    The following is elementary:

\begin{Lemma} \label{opensets} For any pair $(i,j)$, let $U_{i,j}$
represent the open set of $\variety$ where $x^{(0)}_{i,j} \neq
0$, and let $U$ represent the open set of $\variety$ where no
$x^{(0)}_{i,j}$ (for $1 \le i \le m$, $1 \le j \le n$) is zero.
Then $Z_0 = \overline{U_{i,j}} = \overline{U}$, where the bar
represents the closure of the respective sets.
\end{Lemma}
\begin{proof} It is clear that $U_{i,j} \subset Z_0$,
from which it follows that $\overline{U_{i,j}} \subset Z_0$. For
any minimal prime $P$ of $\ideal$ that does not contain
$x^{(0)}_{i,j}$, let $Z(P)$ denote its zero set.  Then $U_{i,j}
\cap Z(P)$ is nonempty, since otherwise, $Z(P) \subset
Z(<x^{(0)}_{i,j}>)$, which would force $x^{(0)}_{i,j} \in P$.
Hence, $U_{i,j} \cap Z(P)$ is dense in $Z(P)$, so
$\overline{U_{i,j}}$ must contain all of $Z(P)$. This shows that
$\overline{U_{i,j}} = Z_0$.  A similar argument shows that
$\overline{U} = Z_0$.
\end{proof}

Now let $Z_1$ represent the subvariety of $\variety$ where all
$x^{(0)}_{i,j}$, for $1 \le i \le m$, $1 \le j \le n$, are zero.
If there are minimal prime ideals of $\ideal$ in $S = F[\vars]$
that contain some $x^{(0)}_{i,j}$ (and hence, by Proposition \ref
{doesntmatterwhichxij} above, contain all $x^{(0)}_{i,j}$ for all
$i,j$), then the zero sets of such prime ideals will clearly be
components of $Z_1$. (These zero sets will also be components of
$\variety$, of course. If no such minimal primes exist, then
$Z_1$ will be contained in $Z_0$.)

Since every minimal prime ideal of $\ideal$ either contains some
(hence all) $x^{(0)}_{i,j}$ or does not contain some (hence any)
$x^{(0)}_{i,j}$, we  have the following:

\begin{Theorem} \label{breakup} The variety $\variety$ is the union of two
subvarieties $Z_0$ and $Z_1$.   The variety $Z_0$ is the closure
of any of the open sets $U_{i,j}$ ($1 \le i \le m$, $1 \le j \le
n$) of $\variety$ where $x^{(0)}_{i,j}$ is nonzero (as also the
closure of the open set $U$ where all $x^{(0)}_{i,j}$ are
nonzero). $Z_0$ is also the union of the components of $\variety$
that correspond to minimal primes of $\ideal$ that do not contain
some (hence do not contain any) $x^{(0)}_{i,j}$. Such components
always exist, and are in one-to-one correspondence with the
components of the variety $\lvariety$. The correspondence
preserves the codimension (in $\affine^{mnk}$ and
$\affine^{(m-1)(n-1)k}$  respectively) of the components. In
fact, $Z_0$ is birational to $\lvariety \times
\affine^{(m+n-1)k}$. The variety $Z_1$ is the subvariety of
$\variety$ where all $x^{(0)}_{i,j}$ are zero, and is isomorphic
to $\Z^{m,n}_{r,k-r} \times \affine^{mn(r-1)} $ when $k
> r$, and isomorphic to $\affine^{mn(k-1)}$ when $k \le r$.

The subvariety $Z_1$ will be wholly contained in $Z_0$ precisely
when there are no minimal primes of $\ideal$ that contain some
(hence all) $x^{(0)}_{i,j}$. If there exist minimal primes of
$\ideal$ that contain some (hence all) $x^{(0)}_{i,j}$, then
these will correspond to some (possibly all) components of $Z_1$.
\end{Theorem}

\begin{proof}
This is just a summary of the discussions in this section.
\end{proof}

\section{The Case of Maximal Minors}
When $r = m$, i.e., when we consider the situation where we set
all maximal minors to zero, we have the following easy result:
\begin{Theorem} \label{maxminorirred}The varieties $\Z^{m,n}_{m,k}$ are all irreducible,
of codimension $k(n-m+1)$.
\end{Theorem}
\begin{proof} We prove the irreducibility by induction on $m$.  If $m=1$, then
the varieties $\Z^{1,n}_{1,k}$ (for any $n$ and $k$) are clearly
irreducible, in fact, $\Z^{1,n}_{1,k}$ is just the origin in
$\affine^{nk}$. So assume that $\Z^{m-1,n-1}_{m-1,k}$ is
irreducible. Then there is only one minimal prime ideal lying
over the ideal $\lideal$ in the ring $T = F[\zvars]$ (see the
statement of Theorem \ref{induct} for notation).  Tracing through
the isomorphism of Theorem \ref{induct} above, the localization of
$\ideal$ at $x^{(0)}_{m,n}$ has only one minimal prime ideal, so
in particular, there is only one minimal prime ideal of $\ideal$
in $S = F[\vars]$ that does not contain $x^{(0)}_{m,n}$. As in the
discussion in Section \ref {reductionsetup} (in particular, see
Theorem \ref{breakup}), this means that subvariety $Z_0$ is
irreducible.  It is now sufficient to show that $Z_1 \subset Z_0$.
We will do this by showing that each point in $Z_1$ is on a line,
all but a finite number of points of which lie inside one of the
open sets $U_{i,j}$. Since $Z_0$ is the closure of any of the
open sets $U_{i,j}$, this will establish that $Z_1 \subset Z_0$.

Let $Q$ be a point in $Z_1$.  If $Q$ is the origin in
$\affine^{mnk}$, then $Q$ lies on the line $\lambda P$ ($\lambda
\in F$) for any $P$ in any $U_{i,j}$.  (Recall that $U_{i,j}$ is
nonempty.)  Since for $\lambda \neq 0$ the point $\lambda P$ is
in $U_{i,j}$, our point $Q$ must lie in the closure of $U_{i,j}$.

Now assume $Q$ is not the origin.  In the representation of $Q$ as
rows $(\u_1(t),\dots,\u_m(t))^T$, with $\u_i(t) = \sum_{l}
\u_i^{(l)}t^l$ (see the notation introduced just before Lemma \ref
{allbasevarszero}), $\u_i^{(0)} = 0$ for $i = 1,\dots,m$. Since
$Q$ is not the origin, some $\u_{i}^{(s)} \neq 0$ for some $i$
with $1 \le i \le m$, and some $s$ with $1 \le s < k$ and with $s$
minimal for this $i$. Write $\vec(t)$ for the vector $\u_{i}^{(s)}
+ \u_{i}^{(s+1)} t + \cdots + \u_{i}^{(k-1)}t^{k-s-1}$ in
$(\trunc)^n$. Consider the point $P(\lambda) = (\u_1(t), \dots,
\u_i(t), \u_{i+1}(t) + \lambda \vec(t),\u_{i+2}(t) \dots,
\u_m(t))^T$ (with the understanding that if $i = m$, then
$P(\lambda) = (\u_1(t), \dots, \u_{m-2}(t),\u_{m-1}(t)+ \lambda
\vec(t), \u_m(t)^T$). Then the $m$-fold wedge product of these
vectors contains two summands:
$\u_1(t)\wedge\dots\wedge\u_i(t)\wedge\u_{i+1}(t)\wedge\dots\wedge
\u_m(t)$, and $\lambda
\u_1(t)\wedge\dots\wedge\u_i(t)\wedge\vec(t)\wedge\dots\wedge
\u_m(t)$ (suitably modified if $i= m$).  The first is zero, since
$Q$ is in $\Z^{m,n}_{m,k}$, and the second is zero since $\u_i(t)
= t^s \vec(t)$.  Thus, the point $P(\lambda)$ is in
$\Z^{m,n}_{m,k}$.  When $\lambda \neq 0$, $P(\lambda)$ is
actually in $U_{i,l}$ for some $l$ (corresponding to any one
coordinate of $\u_{i,s}$ that is nonzero).  Hence, the point $Q =
P(0)$ is in the closure of $U_{i,l}$, which is $Z_0$.

As for the codimension, the codimension of $\I^{m,n}_{m,k}$ is
the codimension of its unique minimal prime ideal.  Tracing
through the localization correspondence of Theorem \ref{induct},
this is the same as the codimension of $\I^{m-1,n-1}_{m-1,k}$.
Continuing this localizing process, we find that the codimension
of   $\I^{m,n}_{m,k}$
is the same as the codimension of $\I^{1,n-m+1}_{1,k}$, which is
clearly $k(n-m+1)$.
\end{proof}

\begin{Remark} Notice how this proof technique breaks down when
considering $r$-fold wedge products  with  $r < m$: if one were
to consider a wedge product that includes the $i$-th row but not
the $(i+1)$-th row (or contains the $m$-th row but not the
$(m-1)$-th row if $i=m$), then the second summand of the wedge
product need not be zero, so the point $P(\lambda)$ need not be
in our variety at all.
\end{Remark}

\subsection{Square Matrices}
When $m=n$, i.e., when our matrices are square, and when we are
still in the situation of maximal minors, we can say considerably
more.   Let us denote the $k$ coefficients of $t$ of the
determinant of our square matrix $X(t)$ by $d_l$, $l=0, \dots,
k-1$. It is easy to determine the structure of the polynomial
expressions $d_l$. The first term $d_0$ is just the determinant
of the matrix $X(0)= ((x^{(0)}_{i,j}))_{1\le i,j\le m}$.  The
remaining terms can be obtained by the following process: Every
monomial appearing in $d_0$ is the form
$$
m_\sigma = x^{(0)}_{\sigma(1),1}x^{(0)}_{\sigma(2),2}\dots
x^{(0)}_{\sigma(m),m}
$$
for some permutation $\sigma$ of $\{1,2,\ldots,m\}$. Given such a
monomial, we define
$$
\mu_s(m_\sigma) = \sum_{k_i\ge 0, \sum k_i = s}
x^{(k_1)}_{\sigma(1),1}x^{(k_2)}_{\sigma(2),2}\dots
x^{(k_m)}_{\sigma(m),m}
$$   We then find
$$
d_s = \sum_{\sigma\in S_n}\mu_s(m_\sigma)sgn(\sigma)
$$

We will prove that the ideals $\sqideal$ are prime  and that
the coordinate rings of the varieties   $\gensqvariety$   are
complete intersection rings and hence Cohen-Macaulay. We will do
so by showing that the $d_l$ form a Groebner basis for $\sqideal$
with respect to a suitable monomial ordering.

We consider the graded reverse lexicographic ordering
(\textsl{grev\-lex}) on the monomials on $S = F[\vars]$ given by
the following scheme: $x^{(k-1)}_{1,1} > x^{(k-1)}_{1,2} > \dots
> x^{(k-1)}_{1,m} > x^{(k-1)}_{2,1} > \dots > x^{(k-1)}_{m,m} >
x^{(k-2)}_{1,1} > x^{(k-2)}_{1,2} > \dots > x^{(k-2)}_{1,m} >
x^{(k-2)}_{2,1} > \dots > x^{(k-2)}_{m,m} > \dots > x^{(1)}_{m,m}
> x^{(0)}_{1,1}
>\dots >x^{(0)}_{1,m} > x^{(0)}_{2,1} > \dots > x^{(0)}_{m,m}$.

(Recall that in the graded reverse lexicographic ordering the
monomials of $S$ are first ordered by the total degree, and for
two monomials $\alpha$ and $\beta$ of the same degree, $\alpha$
is greater than $\beta$ if the rightmost nonzero element in
$\alpha-\beta$ (with $\alpha$ and $\beta$ thought of as elements
of   $\mathbb{Z}^{km^2}$) is negative---see \cite[Chapter 2, \S
2]{CLS}, for instance.)

\begin{Theorem}\label{grob}
Under the grevlex ordering on $S$ described above, the generators
$d_l$, $l=0,1,\ldots,k-1$ of the ideal $\sqideal$   form a
Groebner basis for  $\sqideal$.
\end{Theorem}

\begin{proof}
The grevlex order is designed so as to favor monomials in which
the lower order  variables  do not appear.  It is easy to see
then that the leading monomials ($\lm$) of the various $d_l$ are
as follows:
\begin{eqnarray*} \label{dklt}
\lm(d_0) &=& x^{(0)}_{1,m}x^{(0)}_{2,m-1}\dots x^{(0)}_{m,1}\\
\lm(d_1) &=& x^{(0)}_{1,m-1}x^{(0)}_{2,m-2}\dots x^{(0)}_{m-1,1}x^{(1)}_{m,m}\\
\lm(d_2) &=& x^{(0)}_{1,m-2}x^{(0)}_{2,m-3}\dots x^{(0)}_{m-2,1} x^{(1)}_{m-1,m}x^{(1)}_{m,m-1}\\
\vdots  &=& \vdots\\
\lm(d_{m-1}) &=& x^{(0)}_{1,1} x^{(1)}_{2,m}x^{(1)}_{3,m-1}\dots x^{(1)}_{m,2}\\
\lm(d_m) &=& x^{(1)}_{1,m}x^{(1)}_{2,m-1}\dots x^{(1)}_{m,1}\\
\lm(d_{m+1}) &=& x^{(1)}_{1,m-1}x^{(1)}_{2,m-2}\dots x^{(1)}_{m-1,1}x^{(2)}_{m,m}\\
\vdots  &=& \vdots
\end{eqnarray*}
so in general, we have
\begin{equation} \label{dkltformula}
\lm(d_{\lambda m + \mu}) =
x^{(\lambda)}_{1,m-\mu}x^{(\lambda)}_{2,m-\mu-1}\dots
x^{(\lambda)}_{m-\mu,1}x^{(\lambda+1)}_{m-\mu+1,m}
x^{(\lambda+1)}_{m-\mu+2,m-1} \dots x^{(\lambda+1)}_{m,m-\mu+1},
\end{equation}
where $0\le\mu\le m-1$.

It is clear that the leading monomials  of $d_i$ and $d_j$, for
distinct $i$ and $j$,  are formed of sets of
variables that are disjoint from one another. Hence, the leading
terms of the various $d_k$ are all pairwise relatively prime.  It
follows, e.g. from \cite[Chapter 2, \S9, Prop. 4 and Theorem
3]{CLS}, that the  $d_l$ form a Groebner basis for $\sqideal$.
\end{proof}

We are now ready to prove our   main   result about the ideals
$\sqideal$.

\begin{Theorem} \label{idealsnnl}
The ideals $\sqideal$ are prime ideals, of codimension $k$.  The
coordinate rings of the varieties $\Z^{m,m}_{m,k}$ are
consequently complete intersection rings, and hence
Cohen-Macaulay.
\end{Theorem}
\begin{proof} Since the polynomials $d_l$ form a Groebner basis
for $\sqideal$, and since the lead terms of these $d_l$ in
(\ref{dklt}) are obviously square free, it follows that the
ideals $\sqideal$ are radical.  (This is well known and easy: if
$f^r \in \sqideal$, then $(\lm(f))^r \in
<\lm(d_0),\dots,\lm(d_{k-1})>$, so $\lm(d_i)$ divides
$(\lm(f))^r$ for some $i$.  Since $\lm(d_i)$ is square free, this
means that $\lm(d_i)$ divides $\lm(f)$, so $\lm(d_i) e = \lm(f)$
for some monomial $e$.  Then $f-ed_i$ is also in the radical of
$\sqideal$, and has lower lead monomial.  We proceed thus to find
that $f$ is in $\sqideal$.)  But we have already seen in Theorem
\ref{maxminorirred} above that the radical of $\sqideal$ must be
prime.  It follows that the ideals $\sqideal$ are prime  and that
their codimension is equal to $k$.
Then, by the definition of a complete intersection ring, it
follows that the coordinate ring $S/\sqideal$ of the variety
$\sqvariety$ is a complete intersection ring.  By e.g.
\cite[Chapter 18, Prop. 18.13]{E}, it is Cohen-Macaulay.
\end{proof}

We can obtain more insights into $\sqvariety$ by considering the
quotient ring of $F[\vars]$ by the ideal $\lm(\sqideal)$
generated by the leading monomials of the elements of
$\sqideal$.  By Theorem \ref{grob}, the generators of
$\lm(\sqideal)$ are given by (\ref{dkltformula}) above. As is
well known, the Hilbert function of $\sqideal$ is the same as
that of $\lm(\sqideal)$.  In turn, because $\lm(\sqideal)$ is
generated by the squarefree monomials in  (\ref{dkltformula}), we
may consider the simplicial complex attached to
$F[\vars]/\lm(\sqideal)$, for which   this ring   is the
Stanley-Reisner ring.  (The connection between simplicial
complexes and its associated Stanley-Reisner ring may be found,
for instance, in \cite[Chap. 5 ]{BH}.  Briefly, if $\Delta$ is a
simplicial complex on the vertex set $\{ v_1, \dots, v_n\}$, then
the Stanley-Reisner ring associated to $\Delta$ is the ring
$F[x_1,\dots, x_n]/I_\Delta$, where $I_\Delta$ is generated by
all monomials $x_{i_1}\dots x_{i_s}$, $s\le n$, such that
$x_{i_1}\dots x_{i_s}$ is not a face of $\Delta$.  This
correspondence can be reversed: given an ideal $I$ of
$F[x_1,\dots, x_n]$ such that $I\subset (x_1,\dots,x_n)^2$ and
such that $I$ is generated by squarefree monomials, then we take
$\Delta$ to be the simplicial complex on the vertex set $\{ v_1,
\dots, v_n\}$ whose faces are all those subsets $\{x_{i_1},\dots,
x_{i_s}\}$, $s\le n$, such that $x_{i_1}\dots x_{i_s}$ is not in
$I$.)

We first determine the simplicial complex attached to
$\lm(\sqideal)$.  We start with the following trivial observation:

\begin{Lemma} \label{streisring} The ring $F[\vars]/\lm(\sqideal)$
is isomorphic to $L[x_{km+1},\dots,x_{km^2}]$, where $L$ is
isomorphic to   $F[y_1,\dots,y_{km}]$ modded out by the ideal
generated by $y_1y_2\cdots y_m,$ $y_{m+1}y_{m+2}\cdots y_{2m},$
$\dots,$ $y_{(k-1)m+1}y_{(k-1)m+2}\dots y_{km}$.  (Here,
$y_1,\dots,y_{km}$ are variables.)
\end{Lemma}
\begin{proof}
 Notice from Equality (\ref{dkltformula}) that the leading
monomials of the generators of $\sqideal$ only involve $km$
variables. The form of the leading monomials now gives us the
result above.
\end{proof}

 For any $l \ge 1$, write $C_l$ for the simplicial complex defined by the set of all
subsets of the $l$ element set $\{x_1,\dots,x_l\}$, and write
$S_l$ for the simplicial complex defined by the set of all
subsets of the $l$ element set $\{x_1,\dots,x_l\}$
\textit{except} the full set $\{x_1,\dots,x_l\}$.  (Thus, $S_l$
is the \textsl{skeleton} of the complex $C_l$.)   Let $S_l^k$
be the \textsl{join} of $k$ disjoint copies of $S_l$.  (The join
of two simplicial complexes $\Delta_1$ and $\Delta_2$ on the
disjoint vertex sets $V_1$ and $V_2$ is the complex on the vertex
set $V_1 \cup V_2$ with faces $F_1\cup F_2$, where $F_1$ is  a
face of $\Delta_1$ and $F_2$ is a face of $\Delta_2$.)  We now
have the following essentially trivial result:
\begin{Lemma} \label{simpcomp} The simplicial complex whose
Stanley-Reisner ring is the ring $F[\vars]/\lm(\sqideal)$ is the
complex $S_m^k * C_{k(m^2-m)}$, where $*$ denotes the join.
\end{Lemma}
\begin{proof} This is clear from Lemma \ref{streisring} above. \end{proof}

Recall that the face vector $f(\Delta)$ of a simplicial complex
$\Delta$ of dimension $d-1$ is the $d$-tuple
$(f_0,f_1,\dots,f_{d-1})$, where $f_i$ is the number of $i$
dimensional faces of $\Delta$.  Recall too that a face of
dimension $s$ corresponds to a subset of cardinality $s+1$. Hence,
for the simplicial complex $S_m$, $f_0 = m$, $f_1 =
{{m}\choose{2}}$, $\dots$, $f_{m-2} = {{m}\choose{m-1}}$.  We have
the following:
\begin{Lemma} \label{fvectorsforL}The $f$-vector of the simplicial complex
$S_m^k$ is determined by the coefficients of $x^l$,
$l=1,\dots,k(m-1)$, in the $k$-th power of the polynomial $(1+mx +
{{m}\choose{2}}x^2 + \dots + {{m}\choose{m-1}}x^{m-1})$.
\end{Lemma}
\begin{proof} This is clear from the definition of the join of two
simplicial complexes, and that the fact that the coefficients
$m$, ${{m}\choose{2}}$, etc. in the polynomial above are just the
components of the face vector of $S_m$.  Note that an
$(s-1)$-dimensional face of $S_m^k$ (which corresponds to a
subset of cardinality $s$) arises from a choice of an $(s_1-1)$
dimensional face from the first copy of $S_m$ (corresponding to a
subset of cardinality $s_1$), an $(s_2-1)$ dimensional face from
the second copy of $S_m$ (corresponding to a subset of
cardinality $s_2$), and so on, to a choice of an $(s_k-1)$
dimensional face from the $k$-th copy of $S_m$ (corresponding to a
subset of cardinality $s_k$), where $s_1 + s_2 + \cdots + s_k =
s$.
\end{proof}
In the same vein, we have the following, which gives us the
$f$-vector for the simplicial complex attached to
$F[\vars]/\lm(\sqideal)$:
\begin{Proposition} \label{fvectorsforvariety} Write $b$ for
$k(m^2-m)$. Then, the $f$-vector of the simplicial complex
attached to $F[\vars]/\lm(\sqideal)$ is determined by the
coefficients of $x^l$, $l=1,\dots,k(m^2-1)$, in the polynomial
  $(1+mx + {{m}\choose{2}}x^2 + \dots +
{{m}\choose{m-1}}x^{m-1})^k \cdot (1+bx + {b\choose 2}x^2 + \cdots
+ {b\choose {b-1}}x^{b-1} + x^b)$.
\end{Proposition}
\begin{proof} The proof is similar to the proof of the lemma
above.  The simplicial complex attached to
$F[\vars]/\lm(\sqideal)$ is the join of $S_m^k$ and $C_b$ by
Lemma \ref{simpcomp} above.  The $f$-vectors for $S_m^k$ are
determined by Lemma \ref{fvectorsforL} above. The $f$-vectors of
$C_b$ are easy to determine: the number of $s-1$ dimensional
faces is the number of subsets of a $b$ element set of
cardinality $s$, so it is $b \choose s$.  The rest of the
argument is the same.
\end{proof}

Note immediately that the dimension of the simplicial complex
attached to $F[\vars]/\lm(\sqideal)$ is $k(m-1) + b -1 = k(m^2-1)
- 1$, and that $f_{k(m^2-1)-1} = m^k$.  From the connection
between the face vector of $\Delta$ and the Hilbert polynomial of
the associated Stanley Reisner ring (see \cite[p. 204]{BH} for
instance), and from the fact that the Hilbert function of
$F[\vars]/\lm(\sqideal)$ and $F[\vars]/\sqideal$ are the same, we
have the following:
\begin{Proposition} \label{hp} The
(projective) Hilbert polynomial of $F[\vars]/\sqideal$ is given
by   $H(n) = \sum_{i=0}^{k(m^2-1)-1} f_i{{n-1}\choose i}$,
  where the $f_i$ are the components of the
face vector of the simplicial complex $S_m^k * C_{k(m^2-m)}$
determined by Proposition \ref{fvectorsforvariety}. In
particular, the degree of $\sqvariety$ is $m^k$.
\end{Proposition}
\begin{proof} The degree of $\sqvariety$ is just the coefficient
$f_{k(m^2-1)-1}$ attached to the highest degree term $n-1 \choose
d-1$ in the (projective) Hilbert polynomial when it is expressed
as a linear combination of the polynomials $n-1 \choose i$.
\end{proof}

\begin{Remark} The Hilbert polynomial above also shows the
(projective) dimension of $\sqvariety$ to be $k(m^2-1)-1$, which
is consistent with Theorem \ref{maxminorirred}.
\end{Remark}

\section{Equations for Some Open Sets of $\variety$} \label{APEqns}
In this section, we will derive the key equations that will hold
in certain open sets of our variety  and will enable us to show
that our varieties are reducible when $r < m$ (i.e., we set
submaximal minors to zero).

We start with an elementary and well-known result:
\begin{Lemma} \label{lincomb}
 Let $R$ be a commutative ring, and $R^n$ the free
module of rank $n$. Suppose $\u_1,\dots,\u_{r-1} \in R^n$ are such
that for some $\w_1,\dots,\w_{n-r+1}\in R^n$, the product
$\u_1\wedge\dots\wedge \u_{r-1}\wedge \w_1\wedge\dots\wedge
\w_{n-r+1} \in R^*
$ (i.e., the elements $\u_1,\dots,\u_{r-1}$ can be extended to a
basis of $R^n$).    Here $R^*$ denotes the set of all invertible
elements of $R$.   If $\u_1\wedge\dots\wedge \u_{r-1}\wedge \vec
= 0$ for some element $\vec\in R^n$, then $\vec =
\sum_{i=1}^{r-1}\alpha_i \u_i$ for some $\alpha_i \in R$.
\end{Lemma}
\begin{proof} Since $\u_1\wedge\dots\wedge
\u_{r-1}\wedge \w_1\wedge\dots\wedge \w_{n-r+1}$ is the
determinant of the matrix
$[\u_1,\dots,\u_{r-1},\w_1,\dots,\w_{n-r+1}]$, the hypothesis
shows that this matrix is invertible in $R$. Hence we can find the
unique solution $\x = (\alpha_1,\dots,\alpha_n)^T$ to the equation
$$
[\u_1,\dots,\u_{r-1},\w_1,\dots,\w_{n-r+1}] \x = \vec
$$
by Cramer's rule.  But the assumption $\u_1\wedge\dots\wedge
\u_{r-1}\wedge \vec = 0$   shows that $\alpha_j$ must be zero for
$j = r, \dots, n$, since for such $j$, one of the $\w_i$ will be
replaced by $v$ during the solution.  Hence $\vec =
\sum_{i=1}^{r-1}\alpha_i \u_i$.
\end{proof}

For the rest of this section we will write $R$ for the polynomial
ring $F[t]$, and write $\overline{R}$ for the ring $\trunc$.  Let
$M$ be the free $R$ module of rank $n$, and $\overline{M} =
M\otimes_R \trunc$ be the free $\overline{R}$ module of rank
$n$.   For any $v\in M$, write $\overline{v}$ for the image of
$v$ under the map $M \mapsto \overline{M} = M/t^k M$.

Lemma \ref{lincomb} leads to the following:
\begin{Corollary}\label{lincomzerolevel} Suppose ${\u_1},\dots,{\u_r} \in
\overline{M}$ are such that $\u_1 \wedge \dots\wedge \u_r =0$ and
$\u_{1}^{(0)} \wedge \dots \wedge \u_{r-1}^{(0)} \neq 0$ in
  $\bigwedge^{(r-1)} F^n$   (here, $\u_j = \sum_{l=0}^{k-1} \u_j^{(l)}t^l + t^k M$).
Then, there  are   $\alpha_1,\dots,\alpha_{r-1} \in \overline{R}$
such that $\u_r = \sum_{i=1}^{r-1} \alpha_i \u_i$.
\end{Corollary}
\begin{proof} Since $\u_{1}^{(0)} \wedge \dots \wedge \u_{r-1}^{(0)} \neq 0$
in   $\bigwedge^{(r-1)} F^n$,   there are elements $\w_1,$ $
\dots,$ $\w_{n-r+1}\in F^n$ such that
$\u_{1}^{(0)},\dots,\u_{r-1}^{(0)},\w_1,\dots,\w_{n-r+1}$ form a
basis for   the vector space $F^n$.    In particular,
$\u_{1}^{(0)} \wedge \dots \wedge \u_{r-1}^{(0)} \wedge \w_1
\wedge \dots\wedge \w_{n-r+1} \neq 0$ in $F$.  Thus, the constant
term in  the wedge product $\u_{1} \wedge \dots \wedge \u_{r-1}
\wedge \w_1 \wedge \dots\wedge \w_{n-r+1}$ will be nonzero, so
$\u_{1} \wedge \dots \wedge \u_{r-1} \wedge \w_1 \wedge
\dots\wedge \w_{n-r+1}$ is in $\overline{R}^*$.
The result now follows from Lemma \ref{lincomb}.
\end{proof}

We now come to the main result that generates equations for
certain open sets of our variety:

\begin{Theorem} \label{eqnsforopensets}
Suppose $\u_1, \dots, \u_m \in R^n$ (with $m\le n$) are of degree
at most $k-1$ (that is, each component of each $\u_j$ is a
polynomial of degree at most $k-1$  in $t$), and suppose that
$\overline{\u_{j_1}} \wedge\dots\wedge \overline{\u_{j_r}} = 0$
for all sequences $1 \le j_1 < \dots < j_r \le m$.  Further,
assume that $\u_{1}^{(0)} \wedge\dots\wedge \u_{r-1}^{(0)} \neq
0$ (here, $\u_j = \sum_{l=0}^{k-1} \u_j^{(l)}t^l$).  Then
$$
\u_{j_1} \wedge \dots\wedge \u_{j_r}\wedge \u_{j_{r+1}} \in t^{2k}
\bigwedge^{r+1} M
$$
for all sequences $1 \le j_1 < \dots < j_r < j_{r+1} \le m$.
\end{Theorem}
\begin{proof}
Applying Corollary \ref{lincomzerolevel}  to the elements
$\overline{\u_1}$, $\dots$, $\overline{\u_{r-1}}$,
$\overline{\u_j}$ ($r \le j \le m$), we find
$$
\overline{\u_j} = \sum_{i=1}^{r-1} {\alpha_{j,i}} \overline{\u_i}
$$ for elements $\alpha_{j,i}$ in $\overline{R}$.
If $\alpha_{j,i} = p_{j,i}(t) + t^kR$ for uniquely determined
polynomials $p_{j,i} \in R$ of degree at most $k-1$, define
$$
\vec_j =
  \begin{cases}
    \u_j & {j=1,\dots,r-1}, \\
     \sum_{i=1}^{r-1} {p_{j,i}} {\u_i}& {j = r, r+1,\dots, m}.
  \end{cases}
$$
Then, since the $\vec_j$ depend linearly on the $r-1$ vectors
$\vec_1, \dots, \vec_{r-1}$,   we find, this time in $M$, that
\begin{equation} \label{grandwedge}
\vec_{j_1}\wedge\dots\wedge \vec_{j_r}\wedge \vec_{j_{r+1}} = 0
\end{equation}
 for all sequences $1 \le j_1 < \dots < j_r < j_{r+1} \le m$.

Now, for any $j$, $\u_j$ and $\vec_j$ are equal modulo $t^k$, so
we may write   $\u_j = \vec_j + t^k \y_j$ for suitable $\y_j$.
Then,
$$
\u_{j_1} \wedge \dots\wedge \u_{j_r} \wedge \u_{j_{r+1}} =
(\vec_{j_1} + t^k \y_{j_1}) \wedge \dots \wedge (\vec_{j_{r+1}} +
t^k \y_{j_{r+1}}).
$$
 Expanding the right side, we get the
following: a term $\vec_{j_1} \wedge \dots\wedge \vec_{j_r}
\wedge \vec_{j_{r+1}}$ which is zero by Equality
(\ref{grandwedge}), a sum of terms of the form   $t^k (\pm
\y_{j,i}) \wedge \vec_{j_1} \wedge \dots\wedge
\widehat{\vec_{j_{i}}} \wedge \dots \wedge \vec_{j_{r+1}} $
  (where the hat denotes   the omission
of the term under the hat), and then terms that are clearly in
$t^{2k} \bigwedge^{r+1}M$ or higher. But the product $
 \vec_{j_1} \wedge \dots\wedge \widehat{\vec_{j_{i}}} \wedge
\dots \wedge \vec_{j_{r+1}}$ is already in $t^k \bigwedge^{r}M$,
since on reduction modulo $t^k$, this is just the $r$-fold wedge
product of the various $\overline{\u_j}$.  It follows that
$\u_{j_1} \wedge \dots\wedge \u_{j_r} \wedge \u_{j_{r+1}}$ is in
$t^{2k} \bigwedge^{r+1} M$.  This proves the theorem.
\end{proof}

\begin{Corollary} \label{actualAP}  Assume the $\u_j$ ($1 \le j \le m$) are as in
the preceding theorem.  Then, \begin{equation} \label{APrXr}
\sum_{\stackrel{l_1 + \cdots l_{r+1} = w}{0\le l_j < k}} \u_{j_1}^
{(l_1)}\wedge \dots\wedge \u_{j_r}^{(l_r)} \wedge
\u_{j_{r+1}}^{(l_{r+1})} = 0 \end{equation}
 for each $w$ such that $0 \le w < 2k$, and for
all sequences $1 \le j_1 < \dots < j_r < j_{r+1} \le m$.  In
particular, when $r=2$, \begin{equation} \label{AP2X2}
\sum_{\stackrel{l_1 + l_2 +l_{3} = w}{0\le l_1, l_2, l_3 < k}}
\u_{j_1}^{( l_1)}\wedge \u_{j_2}^{(l_2)}\wedge
\u_{j_{3}}^{(l_{3})} = 0
\end{equation} on the subvariety $Z_0$ of $\Z^{m,n}_{2,k}$, for each $w$
such that $0 \le w < 2k$, and for all sequences $1 \le j_1 < j_2
< j_{3} \le m$.
\end{Corollary}
\begin{proof}    The expressions on the left hand side of
the equalities
(\ref{APrXr}) are merely the coefficients of $t^w$, $0 \le w <
2k$, of $\u_{j_1}\wedge\dots\wedge\u_{j_r}\wedge\u_{j_{r+1}}$, so
by the previous theorem these are zero whenever $\u_{1}^{(0)}
\wedge\dots\wedge \u_{r-1}^{(0)} \neq 0$.  When $r=2$, the
subvariety $Z_0$ is the closure   of   the open set where
$\u_{1}^{(0)} \neq 0$.
\end{proof}

\begin{Remark}\label{moreAPLikeequations} There
is another set of equations that hold on the closure of the open
set of $\variety$ where $\u_{1}^{(0)} \wedge\dots\wedge
\u_{r-1}^{(0)} \neq 0$.  By Corollary \ref{lincomzerolevel}, all
other vectors $\u_i$ can be expressed as an $R$-linear
combination of $\u_1, \dots, \u_{r-1}$.  In particular, all other
vectors $u^{(l)}_i$ can be expressed as an $F$-linear combination
of the $k(r-1)$ vectors $\u_1^{(l)}, \dots, \u_{r-1}^{(l)}$,
$l=0,\dots,k-1$.  It follows that the $(k(r-1)+1)$-fold wedge
product of any of the vectors $\u_i^{(l)}$ must be zero on this
closure, for $i=1,\dots,m$ and for $l=0,\dots,k-1$. This will
   be   trivially true if $(k(r-1)+1) > n$.
\end{Remark}

\section{The Case of $2\times 2$ Submaximal Minors}

In this section, we will completely describe the components of the
variety $\Z^{m,n}_{2,k}$ when $n\ge m \ge 3$ and $k \ge 2$.

For this section, we will write $Y_0$ for our variety
$\Z^{m,n}_{2,k}$, and $X_0$ for the closure $Z_0$ of any of the
open sets $U_{i,j}$ described in Theorem \ref{breakup}.  We will
write $Y_1$ for the subvariety $Z_1$ of Theorem \ref{breakup}
where all $x^{(0)}_{i,j}$ are zero.  We will also write
$\Sigma_0$ for the set $\{0,1,\dots,k-1 \}$.

Now $Y_1$ is isomorphic to $\Z^{m,n}_{2,k-2} \times \affine^{mn} $
when $k
> 2$, and isomorphic to $\affine^{mn}$ when $k = 2$.

In particular, recall from the proof of Lemma
\ref{allbasevarszero} that in the case $k
> 2$, the variety $Y_1$ is really  determined by considering
the generic $m\times n$ matrix with rows $\u_{i}^{(1)} +
\u_{i}^{(2)} t + \cdots + \u_{i}^{(k-2)} t^{k-3}$  ($1 \le i \le
m$) and setting determinants of $2 \times 2$ minors to zero
modulo $t^{k-2}$. We will write $\Sigma_1$ for the  set
$\{1,\dots,k-2 \}$ which indexes the powers of $t$ from each row
of the original matrix which are now governed by
 determinantal equations modulo $t^{k-2}$, and we will write
 $W_1$ for the factor isomorphic to $\Z^{m,n}_{2,k-2}$ determined by
these rows.

Notice  that if $k=3$, the factor $W_1$ is the classical
determinantal variety of $2\times 2$ minors of the generic matrix
$((x^{(1)}_{i,j}))$, and this variety is known to be irreducible
(\cite{BV}). Hence, when $k = 3$, $Y_1$ is just an irreducible
variety cross an affine piece, and is hence irreducible.

\medskip
\centerline{\fbox{ \setlength{\unitlength}{1pt}
\begin{picture}(350,380)
  \put(0,0){\includegraphics{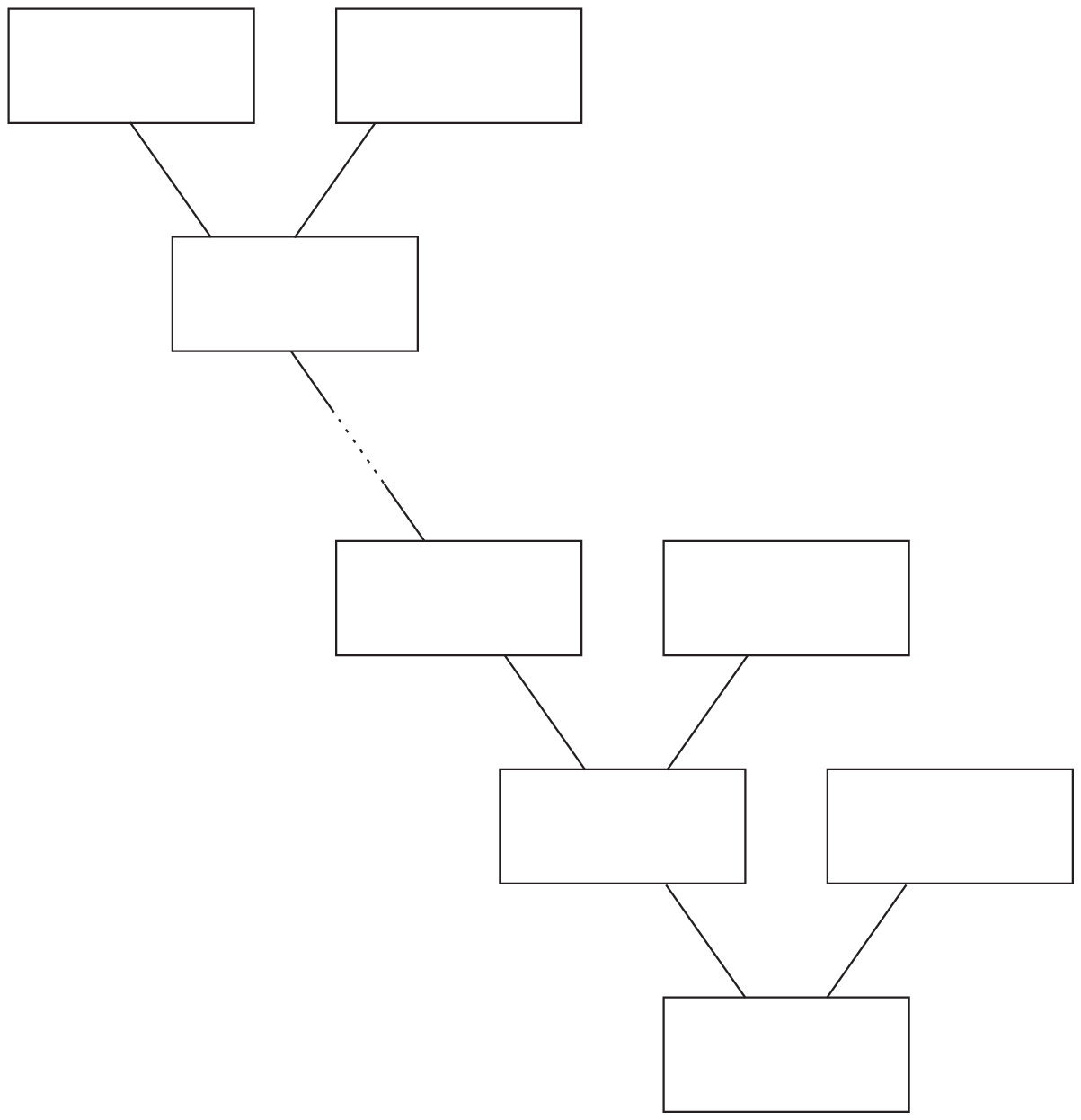}}
  \put(5,335){$\begin{array}{c}Y_L 
  \\ x_{ij}^{(L-1)}=0\end{array}$}
  \put(110,335){$\begin{array}{c}X_{L-1} 
  \\ x_{ij}^{(L-1)}\neq 0\end{array}$}
  \put(60,263){$\begin{array}{c}Y_{L-1} 
  \\ x_{ij}^{(L-2)}=0\end{array}$}
  \put(120,163){$\begin{array}{c}Y_{2} 
  \\ x_{ij}^{(1)}=0\end{array}$}
  \put(227,163){$\begin{array}{c}X_{1} 
  \\ x_{ij}^{(1)}\neq 0\end{array}$}
  \put(170,90){$\begin{array}{c}Y_{1} 
  \\ x_{ij}^{(0)}=0\end{array}$}
  \put(277,90){$\begin{array}{c}X_{0} 
  \\ x_{ij}^{(0)}\neq 0\end{array}$}
  \put(220,15){$\begin{array}{c}Y_{0} \cong \Z_{2,k}^{m,n}
  \end{array}$}
\end{picture}
}}\medskip \centerline{ Figure 1: Irreducible components of
$\Z_{2,k}^{m,n}$}
\medskip

 If $k > 3$, we will write $T_1$
for the subvariety ``$Z_0$'' of $W_1$, i.e., the closure in $W_1$
of the open set where some $x^{(1)}_{i,j}\neq 0$. Also, we will
write $X_1$ for the closure of any of the open sets of $Y_1$
where some $x^{(1)}_{i,j} \neq 0$. It is clear that $X_1$ is just
$T_1 \times \affine^{mn}$.


Write $k = 2L+1$ or $k=2L$ according to whether $k$ is odd or
even. Proceeding thus, we will have subvarieties $Y_s$, $X_s$,
for $s = 0, 1, \dots, L=\lfloor{k/2}\rfloor$. We have the
following:
\begin{itemize}
\item $Y_0$, $X_0$ and
$\Sigma_0 = \{0,1,\dots,k-1\}$ are as already described.
\item (For $0 < s < L$) On $Y_s$, all row vectors
$\u_i^{(0)}, \dots, \u_{i}^{(s-1)}$ are zero, $1 \le i \le m$. The
various rows $\u_{i}^{(s)}, \dots, \u_{i}^{((k-1)-s)}$, $1 \le i
\le m$, are governed by the condition that all $2\times 2$ minors
of the matrix with rows $\u_{i}^{(s)} + \u_{i}^{(s+1)}t + \cdots +
\u_{i}^{((k-1)-s)} t^{k-2s-1} $
 are zero modulo $t^{k-2s}$. We write $\Sigma_s$ for the set
 $\{s,\dots,k-1-s \}$ which indexes  these rows. The
rows $\u_{i}^{(k-s)}, \dots, \u_{i}^{(k-1)}$,  $1 \le i \le m$,
are all free.  Hence, $Y_s \cong W_s \times \affine^{mns}$ $\cong
\Z^{m,n}_{2,k-2s} \times \affine^{mns}$. We   write $T_s$ for the
subvariety ``$Z_0$'' of $W_s$, this is the closure in $W_s$ of
the open set where some $x^{(s)}_{i,j}\neq 0$.
\item (For $0<s < L$) $X_s$   is the closure of any of the open sets of $Y_s$ where
some $x^{(s)}_{i,j} \neq 0$.  It is clear that   $X_s$ is just
$T_s \times \affine^{mns}$.
\item If k = 2L, then $Y_L$ is given by setting all row vectors
$\u_i^{(0)}, \dots, \u_{i}^{(L-1)}$ to zero, $1 \le i \le m$,
with all remaining row vectors $\u_{i}^{(L)}, \dots,
\u_{i}^{(k-1)}$,   $1 \le i \le m$, being free. Thus, $Y_L
\cong\affine^{mnL}$.
\item If k = 2L +1, then on $Y_L$,  all row vectors
$\u_i^{(0)}, \dots, \u_{i}^{(L-1)}$ are zero, $1 \le i \le m$.
The row vectors   $\u_i^{(L)}$,   $1 \le i \le m$, satisfy the
classical determinantal equations for the $2\times 2$ minors of
the generic $m\times n$ matrix $((x^{(L)}_{i,j}))$.   $W_L$ will
denote the classical determinantal variety determined by the
$\u_i^{(L)}$.   The various row vectors $\u_{i}^{(L+1)}, \dots,
\u_{i}^{(k-1)}$, $1 \le i \le m$, are all free. Thus, $Y_L\cong
W_L \times \affine^{mnL}$, and since $W_L$ is a classical
determinantal variety, $Y_L$ is itself irreducible (\cite{BV}).
\end{itemize}

Our result is the following:
\begin{Theorem} \label{2X2decompose}
The variety $\Z^{m,n}_{2,k}$ ($n \ge m \ge 3$, $k\ge 2$) is
reducible.  Its components are the subvarieties $X_0$, $X_1$,
$\dots$, $X_{L-1}$, and $Y_L$ described above.  The components
$X_s$, $s=0,1,\dots, L-1$, have codimension $(m-1)(n-1)(k-2s) +
mns$.  If $k =2L$, $Y_L$ has codimension $mnL$, while if
$k=2L+1$, $Y_L$ has codimension $ (m-1)(n-1) +mnL$.
\end{Theorem}

\begin{proof} The fact that the various $X_s$ and $Y_L$ are irreducible
and have the stated  codimension follows easily from the
descriptions of the various $X_s$ and $Y_L$ above and from
Theorem \ref{induct}. We have seen that for $0 \le s < L$, $X_s
\cong T_s \times \affine^{mns}$, and that $X_s$ sits in the
portion of $\affine^{mn(k-2s)}$  determined by setting all rows
$\u^{(l)}_{i}$ to zero, $l = 0,1,\dots, s-1$ (when $s=0$ this
condition is vacuous). Recall, too, that $T_s$ is the closure of
the open set of $ \Z^{m,n}_{2,k-2s}$ where some $x^{(s)}_{i,j}
\neq 0$. By Theorem \ref{induct}, we have a one-to-one
correspondence between the components of $T_s$ and the components
of $\Z^{m-1,n-1}_{1,k-2s}$, a correspondence which preserves
codimension in the respective spaces $\affine^{mn(k-2s)}$ and
$\affine^{(m-1)(n-1)(k-2s)}$. Since $\Z^{m-1,n-1}_{1,k-2s}$ is
clearly irreducible of codimension $(m-1)(n-1)(k-2s)$ (it is the
origin in $\affine^{(m-1)(n-1)(k-2s)}$), we find that each $X_s$
is irreducible. It follows too that $X_s$ has codimension
$(m-1)(n-1)(k-2s) + mns$ in $\affine^{mnk}$, where the extra
summand $mns$ accounts for the rows $\u_i^{(0)}\dots
\u_i^{(s-1)}$ set to zero.

As for $Y_L$, we have already observed in the discussion before
this theorem that it is irreducible.  In the case $k = 2L$, $Y_L$
is just $\affine^{mnL}$.  It follows that $Y_L$ has codimension
$mnL$ (corresponding to the rows $\u_i^{(0)}\dots \u_i^{(L-1)}$
set to zero).  If $k=2L+1$, then  the codimension of $Y_L$ is
$mnL$ plus the codimension of the variety $W_L$.  But this is
known to be $(m-1)(n-1)$ (see \cite{BV} for instance, this can
also be derived from Theorem \ref{induct}).

We will now prove  that the components of $\Z^{m,n}_{2,k}$ are as
described. It is easily seen from the codimension formulas above
that except when $(m,n) = (3,3)$ or $(m,n)=(3,4)$,  the
codimension decreases as a function of $s$.  This shows that if
$s' > s$, then $X_s'$ (or $Y_L$) cannot be contained in $X_{s}$.
Since the reverse containment is ruled out as $x_{i,j}^{(s)} \neq
0$ on $X_s$, we  find that the components of $\Z^{m,n}_{2,k}$ are
indeed as described, except in the two special cases.

To take care of these two special cases, we will use results from
Section \ref{APEqns}.  (The proof works for all $(m,n)$ pairs
actually.) Using reverse induction on $s$, we will show that at
the $s$-th stage, $s = L, L-1, \dots, 0$, the components of $Y_s$
are $X_s, X_{s+1}$, $\dots$, $X_{L-1}$, and $Y_L$.  We have
already observed that $Y_L$ is irreducible, so assume that
$s<L$.  Note that $Y_s$ is the union of $X_s$ and $Y_{s+1}$, and
by induction, $Y_{s+1}$ will have components $ X_{s+1}$, $\dots$,
$X_{L-1}$, and $Y_L$. We will prove that none of these
subvarieties $ X_{s+1}$, $\dots$, $X_{L-1}$, and $Y_L$ can be
contained in $X_s$.  The reverse containment is once again ruled
out, and we will indeed find that $Y_s$ has components $X_s,
X_{s+1}$, $\dots$, $X_{L-1}$, and $Y_L$.

We will first show that no $X_{s+\alpha}$ ($\alpha = 1, \dots,
L-s-1$) can be contained in $X_s$.  For, assume to the contrary.
Recall that $X_{s+\alpha}$ decomposes as $T_{s+\alpha} \times
\affine^{mn(s+\alpha)}$, where the factor $T_{s+\alpha}$ comes
from the various entries of the rows $\u_i^{(l)}$, $1\le i \le
m$, $l = s+\alpha, \dots, k-1-(s+\alpha)$, and the factor
$\affine^{mn(s+\alpha)}$ comes from the various entries of the
rows $\u_i^{(l)}$, $1\le i \le m$, $l = k-(s+\alpha), \dots,k-1$.
Recall too that $T_{s+\alpha}$ is the closure in $W_{s+\alpha}$
where some $\u_i^{(s+\alpha)} \neq 0$, where $W_{s+\alpha}$ has
the description given earlier. The following point $P$ is
therefore in $X_{s+\alpha}$: $\u_1^{(s+\alpha)} = (1,0,0,\dots)$,
$\u_2^{(k-(s+\alpha))} = (0,1,0,\dots)$, $\u_3^{(k-1-s)} =
(0,0,1,0,\dots)$, and all other rows of all possible degrees
zero. (The nonzero coordinates coming from the row
$\u_1^{(s+\alpha)}$ belong to $T_s$ while those coming from
$\u_2^{(k-(s+\alpha))}$ and $\u_3^{(k-1-s)}$ belong to the other
factor $\affine^{mn(s+\alpha)}$.)

Now $X_s$ decomposes as $T_s \times \affine^{mns}$, where the
factor $T_s$ comes from the various entries of  $\u_i{(l)}$, $1\le
i \le m$, $l = s, \dots, k-1-s$.  Since $P\in X_s$ by assumption,
an examination of its coordinates show that $P$ is in the
subvariety $T_s \times \mathbf{O}$, where we have written
$\mathbf{O}$ for the origin in $\affine^{mns}$.  Thus, the
coordinates of $P$ coming from the rows $\u_i{(l)}$, $1\le i \le
m$, $l = s, \dots, k-1-s$ must satisfy Equations (\ref{AP2X2}) in
Corollary \ref{actualAP} which   hold on $T_s$.  In particular,
the specific equation in (\ref{AP2X2}) that holds for the
coefficient of $t^{2(k-2s)-1}$, on specializing to the rows
$\u_1$, $\u_2$, and $\u_3$, reads

\begin{equation} \label{forXsgen}
\sum_{\stackrel{a+b+c = 2(k-2s)-1}{0 \le a,b,c < k-2s}} \u_{1}^{(
s+a)} \wedge \u_{2}^{(s+b)} \wedge \u_{3}^{(s+c)} = 0
\end{equation}

Examining the coordinates of $P$, and recognizing that
$\u_2^{(k-(s+\alpha))} = \u_2^{(s+(k-2s-\alpha))}$ and
$\u_3^{(k-1-s)} = \u_3^{(s+ (k-1-2s))}$ we find that the only
wedge product that is nonzero in this equation is
$\u_1^{(s+\alpha)} \wedge \u_2^{(k-(s+\alpha))} \wedge
\u_3^{(k-1-s)}$.  But by our choice of these rows, this wedge
product is clearly nonzero.  This shows that $X_{s+\alpha}$ is
not contained in $X_s$.

To show that $Y_L$ is not contained in $X_s$, consider first the
case $k = 2L + 1$.  Then $Y_L \cong \Z^{m,n}_{2,1} \times
\affine^{mnL}$, where the factor $\Z^{m,n}_{2,1}$ comes from the
entries of $\u_i{(L)}$, $1\le i \le m$, and the factor
$\affine^{mnL}$ comes from the entries of $\u_i{(l)}$, $1\le i
\le m$, $l = L+1,\dots,2L$.  Choose $P$ to be the point with
$\u_1^{(L)} = (1,0,0,\dots)$, $\u_2^{(L+1)} = (0,1,0,\dots)$,
$\u_3^{(k-1-s)} = (0,0,1,0,\dots)$, and all other rows of all
possible degrees zero. (The nonzero coordinates coming from the
row $\u_1^{(L)}$ belong to $\Z^{m,n}_{2,1}$ while those coming
from $\u_2^{(L+1)}$ and $\u_3^{(k-1-s)}$ belong to the other
factor $\affine^{mnL}$.) Exactly as before, this point $P$ is in
$Y_L$, but we find that  $P$ does not satisfy Equation
(\ref{AP2X2}) for $T_s$.

When $k=2L$, we take $P$ to be the point with $\u_1^{(L)} =
(1,0,0,\dots)$, $\u_2^{(L)} = (0,1,0,\dots)$, $\u_3^{(k-1-s)} =
(0,0,1,0,\dots)$, and all other rows of all possible degrees
zero.  Once again, this point $P$ is in $Y_L$ but does not satisfy
Equation (\ref{AP2X2}) for $T_s$.  This completes the proof.

\end{proof}

We get the following corollary from this:
\begin{Corollary} \label{2X2codimensions} The variety
$\Z^{m,n}_{2,k}$ has $1 + \lfloor{k/2}\rfloor$ components. The
codimension of $\Z^{m,n}_{2,k}$ (in $\affine^{mnk}$) is
$(m-1)(n-1) + mn\lfloor{k/2}\rfloor$ if $k$ is odd and   the
codimension is $mn\lfloor{k/2\rfloor}$ if $k$ is even, except in
the case where $(m,n) = (3,3)$ or $(m,n)=(3,4)$. In these two
special cases, $\Z^{3,3}_{2,k}$ has codimension $4k$, while
$\Z^{3,4}_{2,k}$ has codimension $6k$.
\end{Corollary}
\begin{proof}  By Theorem \ref{2X2decompose}, $\Z^{m,n}_{2,k}$ clearly has
$1 + \lfloor{k/2}\rfloor$ components.  These have codimension
$(m-1)(n-1)(k-2s) + mns$, for $s=0,1,\dots, L
=mn\lfloor{k/2}\rfloor$ (the case $s=L$ corresponds to the
component $Y_L$, but is also covered by this formula).  This is
linear in $s$, and as already observed, is decreasing in $s$
except when $(m,n) = (3,3)$ or $(m,n)=(3,4)$.
It follows that except for these two cases, the component $Y_L$,
has the least codimension.  This yields the formula for the
codimension of $\Z^{m,n}_{2,k}$ in the general case. In the two
special cases, the component $X_0$ must have least codimension.
Hence, the codimension of $\Z^{m,n}_{2,k}$ in these cases is
given by the codimension of $X_0$, which is $(m-1)(n-1)k$.

\end{proof}

\begin{Remark} \label{34equalcd} Note that when $(m,n)=(3,4)$,
the codimension of the components is constant in $s$.
 Hence, in this case, all components have the
same dimension.
\end{Remark}

\section{Submaximal Minors: The General Situation}

In this section, we will use Theorem \ref{2X2decompose} as a
building block to derive results about $\variety$ in general,
when $r < m$.  The first one is easy:

\begin{Theorem} \label{mincomponentssubmaximalcase} The variety
$\variety$ in the submaximal case ($r < m$) has at least $1 +
\lfloor{k/2}\rfloor$ components.
\end{Theorem}

\begin{proof} By Theorem \ref{breakup}, $\variety$ has at least as
many components as its subvariety $Z_0$, and the components of
this subvariety are in one-to-one correspondence with those of
$\lvariety$. Proceeding thus, $\variety$ has at least as many
components as $\Z^{m-r+2,n-r+2}_{2,k}$, and this last variety has
$1 + \lfloor{k/2}\rfloor$ components.
\end{proof}

It is quite clear that the components of $\Z^{m,n}_{2,k}$
intersect one another.  For instance, the origin is in $Y_L$ and
is also in each of the $X_s$. (For, as in the proof of Theorem
\ref{maxminorirred}, the line between the origin and any point
  $P\in T_s$   with some
$x^{(s)}_{i,j} \neq 0$ has to lie in $X_s$.) But, as well, it is
quite easy to find lots of other points of intersection between
the various components. For example, if $k\ge 4$, then there are
at least three components, $X_0$, $X_1$ and $Y_L$. The rows
$\u_{i}^{(k-1)}$ are free on both $X_1$ and $Y_L$, and it is then
easy to see that the point with zeros in all rows except in the
rows $\u_{i}^{(k-1)}$ is in both $X_1$ and $Y_L$.  If $k=2$, then
the components are $X_0$ and $Y_1$.  For any nonzero $\lambda\in
F$ and for any point    $P \in X_0$ with some $\u_i^{(0)} \neq
0$, the point with coordinates $\lambda\u_i^{(0)}, \u_{i}^{(1)}$,
$i = 1, \dots, m$, also satisfies the equations of
$\Z^{m,n}_{2,2}$, so when
$\lambda$   
equals zero, we get a point in  $X_0 \cap Y_1$.  Similarly, if
$k=3$, then the point with coordinates $\lambda^2\u_i^{(0)},
\lambda\u_{i}^{(1)}, \u_{i}^{(2)}$, $i = 1, \dots, m$, also
satisfies the equations of $\Z^{m,n}_{2,3}$, so once again, if we
start with a point in $X_0$ with some $\u_i^{(0)} \neq 0$ and let
$\lambda$ equal zero, we get a point in $X_0\cap Y_1$.

From this, we get the following trivially:
\begin{Theorem} \label{submaximalNoCMNoNormal} The varieties
$\variety$ when $r < m$ are not normal.  They are also not
Cohen-Macaulay, except possibly in the case where $(m,n) =
(1+r,2+r)$.
\end{Theorem}
\begin{proof} In the case of $2\times 2$ minors, the fact that the varieties
are not normal follows from the fact that we have explicitly
found components
 that intersect nontrivially, while the fact
that they are not Cohen-Macaulay except possibly when
$(m,n)=(3,4)$ follows from the fact that components with
different dimensions intersect nontrivially (see Remark
\ref{34equalcd}). In the general case, we repeatedly invoke the
birational isomorphism of Theorem \ref{induct} between the
subvariety $Z_0$ (which is a union of some of the   components of
$\variety$) and the variety  $\lvariety \times
\affine^{k(m+n-1)}$, and then  reduce to the case
$\Z^{m-r+2,n-r+2}_{2,k}$.  Note that the image
  of the birational isomorphism is the open
 subset of $\lvariety \times
\affine^{k(m+n-1)}$  where the free variable $x^{(0)}_{m,n}\neq
0$ (see Remark \ref{birational}), so any intersection between
components of $\lvariety$ of dimensions $d_1$ and $d_2$ will
indeed manifest itself as an intersection between components  of
$\variety$ (in fact of the subvariety $Z_0$ of $\variety$) of
dimensions  $d_1 +k(m+n-1)$ and $d_2 +k(m+n-1)$.
\end{proof}

For $r \ge 3$ (and $r < m$), it is somewhat more difficult (and
cumbersome) to determine explicitly the components of
$\variety$.  Recall that our variety $\variety$ decomposes into
$Z_0$ and $Z_1$ as in Theorem \ref{breakup}.  The components of
$Z_0$ are in one-to-one correspondence with those of $\lvariety$,
while the components of $Z_1$ are determined by those of
$\Z^{m,n}_{r,k-r} $ when $k > r$, and of course, when $k\le r$,
$Z_1$ is isomorphic to $\affine^{mn(k-1)}$ (by Lemma
\ref{allbasevarszero}).  Now, even if we are able to determine
inductively the components of $\lvariety$ and those of
$\Z^{m,n}_{r,k-r} $ (when $k > r$), it is quite difficult to
determine how these various components sit inside $\variety$.
Specifically, it is quite difficult to determine if any of the
components of $Z_1$ live inside any of the components of $Z_0$.
(It is clear that no component of $Z_0$ can live inside any
component of $Z_1$, since every component of $Z_0$ contains
points with $x^{(0)}_{i,j} \neq 0$.) But there is one special
situation where one can indeed answer this question explicitly,
and this is in the case where $k < r$. We have the following:
\begin{Proposition} \label{k<rreduction} In the case where $k < r$, the subvariety $Z_1$ of
$\variety$ is contained in $Z_0$.  The components of $\variety$
and their codimensions in $\affine^{mnk}$ in this case are hence
determined completely by the components of $\lvariety$ and their
codimensions in $\affine^{(m-1)(n-1)k}$.
\end{Proposition}

\begin{proof} Consider the subvariety $V$ of $\variety$ defined by
setting all $2\times 2$ minors of degree zero to zero, i.e.,
defined by setting all $\u_i^{(0)}\wedge\u_j^{(0)} = 0$ for all $1
\le i < j \le m$.  The equations for $V$ are thus the standard
equations for $\variety$ along with all $2\times 2$ minors of
degree zero. It is clear that every point on $Z_1$ satisfies
these equations, so $Z_1 \subset V$.  Since $r \ge (k-1) + 2$,
every $r$-fold wedge product of vectors $\u_i^{(l)}$ of total
degree at most $k-1$ must contain at least two factors of degree
zero. It follows that $\ideal$ is already contained in the ideal
generated by all $2\times 2$ minors of degree zero, which is an
ideal that is known classically to be prime, \cite{BV}. Hence,
$V$ is an irreducible variety, isomorphic to $\Z^{m,n}_{2,1}.$ The
components of $\variety$ come from either $Z_1$ or $Z_0$. Since
$V$ cannot be contained wholly in any component of $Z_1$ (as there
are clearly points on $V$ where not all $x^{(0)}_{i,j}$ are
zero), we find  $V\subset Z_0$.  It follows that $Z_1 \subset
Z_0$.  Thus the components of $\variety$ all come from $Z_0$, and
Theorem \ref{breakup} now finishes the proof.
\end{proof}

We use the result above to determine the components of the
tangent bundle to the classical determinantal varieties in the
case of submaximal minors:
\begin{Corollary} \label{tanbundlecomp} When $k=2$ (i.e., when we
consider the tangent bundle to $\Z^{m,n}_{r,1}$), and when $r <
m$, $\Z^{m,n}_{r,2}$ has exactly two components.  One of them is
the closure of any of the open sets $U_{[i_1,\dots, i_{r-1} |
j_1,\dots, j_{r-1}]}$ of $\variety$, where  the
$(r-1)\times(r-1)$ minor of degree zero determined by rows
$i_1,\dots, i_{r-1} $ and columns $ j_1,\dots, j_{r-1}$ is
nonzero, and hence also of their union.   This component has
codimension $2(m-r+1)(n-r+1)$. The other is the subvariety
defined by setting all $(r-1)\times (r-1)$ minors of degree zero
to zero, and has codimension $(m-r+2)(n-r+2)$.
\end{Corollary}
\begin{proof}
The number of components and their codimension comes from
repeated applications of Proposition \ref{k<rreduction} above. By
the proposition, the components of $\Z^{m,n}_{r,2}$ are in
one-to-one correspondence (via the birational isomorphism of
Theorem \ref{induct}) to the components of
$\Z^{m-1,n-1}_{r-1,2}$.  If $r-1=2$, then the components of
$\Z^{m-1,n-1}_{r-1,2}$  and their codimensions  are described by
Theorem \ref{2X2decompose}. Otherwise, we repeat the process,
until we come to $\Z^{m-r+2,n-r+2}_{2,2}$.

It remains to establish the description of the components.  We
proceed by induction on $r$. When $r=2$, this is precisely the
content of Theorem \ref{2X2decompose}, as also of Lemma
\ref{opensets}. Note that   if $P_0$ is the prime corresponding
to $Z_0$, then $P_0$ does not contain any $x^{(0)}_{i,j}$, while
if $P_1$ is the prime corresponding to $Z_1$, then $P_1$ is the
minimal prime of $\I^{m,n}_{2,2}$ that contains all
$x^{(0)}_{i,j}$ (Lemma \ref{doesntmatterwhichxij}). Now assume
that the theorem is true for the variety $\Z^{m-1,n-1}_{r-1,2}$.
We will invoke the notation of Theorem \ref{induct}.   The
assumption shows that there are exactly two minimal primes of
$\I^{m-1,n-1}_{r-1,2}$ in $T$. One of them, call it $Q_0$, does
not contain any $(r-2)\times(r-2)$ minors of the degree zero
matrix $((z^{(0)}_{i,j}))$ ($1 \le i \le m-1, 1 \le j \le n-1$),
since otherwise, the component $Z(Q_0)$ corresponding to $Q_0$
cannot contain the open set where this minor is nonzero. The
other, call it $Q_1$ is the unique minimal prime over the ideal
of $T$ generated by $\I^{m-1,n-1}_{r-1,2}$ and the various
$(r-2)\times(r-2)$ minors of the degree zero matrix
$((z^{(0)}_{i,j}))$.   Let $\widetilde{Q}_0$ and
$\widetilde{Q}_1$ in $S[({x^{(0)}_{m,n}})^{-1}]$ be the images
(respectively) of $Q_0 T' [({x^{(0)}_{m,n}})^{-1}]$ and $Q_1 T'
[({x^{(0)}_{m,n}})^{-1}]$   under the inverse map $\tilde{f}$.
Notice that the various $(r-2)\times(r-2)$ minors of the degree
zero matrix $((z^{(0)}_{i,j}))$ are just the generators of the
ideal $\I^{m-1,n-1}_{r-2,1}$.  By the proof of Theorem
\ref{induct}, these generators go to the generators of
$\I^{m,n}_{r-1,1}[({x^{(0)}_{m,n}})^{-1}]$ under the inverse map
$\tilde{f}$.  Thus,   $\widetilde{Q}_1$ contains all generators
of $\I^{m,n}_{r-1,1}[({x^{(0)}_{m,n}})^{-1}]$.  On the other
hand, the standard generators for
$\I^{m,n}_{r-1,1}[({x^{(0)}_{m,n}})^{-1}]$   are those that come
from $\I^{m,n}_{r,1}$: these are just the various $(r-1)\times
(r-1)$ minors  of our $m\times n$ degree zero matrix
 $X(0)=((x^{(0)}_{i,j}))$.   Thus, the
pullback of   $\widetilde{Q}_1$,  call it $P_1$, contains all
$(r-1)\times (r-1)$ minors of degree zero of $X(t)$. Moreover,
since  $\widetilde{Q}_1$ is the unique minimal ideal of
$S[({x^{(0)}_{m,n}})^{-1}]$   lying over the ideal generated by
$\I^{m,n}_{r,1}$ and the $(r-1)\times (r-1)$ minors of degree
zero of $X(t)$, $P_1$, being its pullback, is the unique minimal
prime of $\I^{m,n}_{r,1}$  lying over the ideal generated by
$\I^{m,n}_{r,1}$ and the $(r-1)\times (r-1)$ minors of degree
zero of $X(t)$.  Thus, the component corresponding to $P_1$ is
the subvariety of $\Z^{m,n}_{r,2}$ obtained by setting all
$(r-1)\times (r-1)$ minors of degree zero to zero.

As for the description of $P_0$, it is clear that $P_0$ cannot
contain any   $(r-1)\times (r-1)$ minor   $[i_1,\dots, i_{r-1} |
j_1,\dots, j_{r-1}]$ of degree zero, since this minor is already
in $P_1$ and would hence be   zero on all of $\Z^{m,n}_{r,2}$ if
it were also in $P_0$, a contradiction.  Recasting this in the
language of open sets, for any $(r-1)\times (r-1)$ minor of
degree zero of $X(t)$ indexed by rows $i_1,\dots, i_{r-1}$ and
columns $j_1, \dots, j_{r-1}$, we must have $U_{[i_1,\dots,
i_{r-1} | j_1,\dots, j_{r-1}]} \subset Z(P_0)$, since
$\Z^{m,n}_{r,2} = Z(P_0)\cup Z(P_1)$ and since the minor
$[i_1,\dots, i_{r-1} | j_1,\dots, j_{r-1}]$ is zero on $Z(P_1)$.
Because $Z(P_0)$ is irreducible and $U_{[i_1,\dots, i_{r-1} |
j_1,\dots, j_{r-1}]}$ is clearly nonempty, the closure of
$U_{[i_1,\dots, i_{r-1} | j_1,\dots, j_{r-1}]}$ must be all of
$Z(P_0)$.  This must then trivially be true of the  union of all
these open sets indexed by these minors. (As well, this is true
of their intersection, as the intersection is also nonempty).
\end{proof}

\begin{Remark}
Except when $(m,n) = (1+r,1+r)$ or $(m,n) = (1+r,2+r)$, the
corollary can be obtained very easily without recourse to the
machinery of this paper. Write $U$ for the union of the open sets
where some $(r-1)\times (r-1)$ degree zero minor of $X(0)$ is
nonzero. Note that the portion of the classical degree zero
variety $\Z^{m,n}_{r,1}$ where some $(r-1)\times (r-1)$ minor is
nonzero is precisely the set of smooth points of
$\Z^{m,n}_{r,1}$. The variety $\Z^{m,n}_{r,2}$ is the union of
two subvarieties: one, call it $X$, is the closure of $U$, and
the other, call it $Y$ is the subvariety where all $(r-1)\times
(r-1)$ degree zero   minors   of $X(0)$ are zero.  It is easy to
see that $U$ is irreducible of the stated codimension, since the
fibers over any point of $\Z^{m,n}_{r,1}$ where some $(r-1)\times
(r-1)$ is nonzero are all linear spaces of the same dimension.
Hence $X$ is irreducible of the stated codimension. It is also
easy to see that the Jacobian matrix defining tangent spaces to
classical variety $\Z^{m,n}_{r,1}$ is zero when all $(r-1)\times
(r-1)$ minors are zero, so indeed, the tangent spaces at such
points are simply copies of $\affine^{mn}$.  Since the base space
$\Z^{m,n}_{r-1,1}$ is irreducible, $Y$ is irreducible as well,
and it has the stated codimension. It is clear that $X$ cannot be
contained in $Y$ as $U$ is nonempty. Except for the given
exceptional values of $(m,n)$, the dimension of $Y$ is greater
than that of $X$, so $Y$ cannot be contained in $X$ as well.  It
follows that $X$ and $Y$ are precisely the components of
$\Z^{m,n}_{r,2}$.
\end{Remark}

We end this section with an inductive scheme for computing the
codimension of $\variety$ in the case $r < m$. The induction is
based on $r$, and we will assume that for all $r'$ with $2 \le r'
< r$ and for all $m$, $n$ with $r' < m \le n$, and for all $k \ge
2 $, we know the codimension of $\Z^{m,n}_{r',k}$. (The starting
point for the induction is Theorem \ref{2X2decompose}, and the
ideas here parallel the codimension computations of Theorem
\ref{2X2decompose}.)

Write $k = \lambda r + \mu$, for $\lambda \ge 1$, and $0 \le \mu
< k$.  (When $k < r$, we already know that the components, and
their codimensions, are determined by those of $\lvariety$,
thanks to Proposition \ref{k<rreduction} above.) We now have the
following sequence of subvarieties:

\begin{itemize}
\item We will write $Y_0$ for our variety $\variety$, and $X_0$
for its subvariety $Z_0$.  Thus, $X_0$ is birational to
$\lvariety \times \affine^{k(m+n-1)}$.  Write $c_0$ for the
codimension of $\lvariety$ in $\affine^{(m-1)(n-1)k}$. Then $X_0$
also has codimension $c_0$.  We will assume that $c_0$ is known
by induction.
\item We will write $Y_1$ for the subvariety $Z_1$ of
$\variety$---this is obtained by setting all $x^{(0)}_{i,j}$ to
zero.  $Y_1$ is isomorphic to $\Z^{m,n}_{r,k-r} \times
\affine^{mn(r-1)}$. We will write $X_1$ for the subvariety
``$Z_0$'' of $Y_1$.  It is birational to $\Z^{m-1,n-1}_{r-1,k-r}
\times \affine^{m+n-1}$.   Write $c_1$ for the codimension of
$\Z^{m-1,n-1}_{r-1,k-r}$ in $\affine^{(m-1)(n-1)(k-r)}$.  Then
the codimension of $X_1$ is $c_1 + mn$, where the extra
codimension $mn$ comes from the fact that $X_1$ sits in the
portion of $\affine^{mnk}$ where all $x^{(0)}_{i,j}$ are zero. We
will assume that $c_1$ is known.
\item Proceeding thus, let $Y_s$ ($s=1,\dots, \lambda-1$) be the subvariety of $Y_{s-1}$
where all $x^{(s-1)}_{i,j}$ are zero, and, let $X_s$ be the
subvariety ``$Z_0$'' of $Y_s$.  Then $X_s$ has codimension $c_s +
smn$, where $c_s$ is the codimension of $\Z^{m-1,n-1}_{r-1,k-rs}$
in $\affine^{(m-1)(n-1)(k-rs)}$.  We will assume that $c_s$ is
known.
\item If $\mu=0$, i.e., if $k = \lambda r$, then $Y_\lambda$,
the subvariety of $Y_{\lambda-1}$ where all
$x^{(\lambda-1)}_{i,j}$   are zero,
 is
already an affine space of codimension $mn\lambda$ in
$\affine^{mnk}$. For convenience we will take   $c_\lambda = 0$
in this case, so the codimension of $Y_\lambda$ may be written
for this case as $c_\lambda + \lambda mn$.
\item If $\mu > 0$, then   $Y_{\lambda}$  is isomorphic
to $\Z^{m,n}_{r,\mu}\times \affine^{\lambda(r-1)mn}$. Since $\mu
< r$, we can reduce its codimension computations to that of
$\Z^{m-1,n-1}_{r-1,\mu}$ by Proposition \ref{k<rreduction}, so we
will assume that $c_{\lambda}$, the codimension of
$\Z^{m,n}_{r,\mu}$ is known. It follows that the codimension of
$Y_{\lambda}$ is $c_{\lambda} + \lambda mn$.
 ($X_\lambda$ in this case will equal
$Y_\lambda$ by Proposition \ref{k<rreduction}.)
\end{itemize}

Our result is the following:
\begin{Theorem} \label{codimingeneral} The codimension of
$\variety$ in the case $r < m$ is the minimum of the numbers $c_s
+ s mn$, $s = 0, 1, \dots, \lambda$.
\end{Theorem}
\begin{proof}
We will show that the codimension of the subvariety $Y_s$,
$s=0,1,\dots,\lambda$ is the minimum of $c_{s'} + s' mn$,
$s'=s,s+1,\dots,\lambda$.
 When $s=\lambda$ the result
is clear. For   $s < \lambda$, the codimension of $Y_{s+1}$, by
reverse induction, is the minimum of $c_{s'} + s'mn$, $s' = s+1,
\dots, \lambda$.
 For such $s$, the
codimension of $Y_s$ is the minimum of the codimension of
$Y_{s+1}$ and $X_s$, since $Y_{s}$ is the union of $Y_{s+1}$ and
$X_s$.    Putting this together with our inductive result, the
codimension of $Y_s$,  is the minimum of $c_{s'} + s'mn$, $s' = s,
s+1, \dots, \lambda$.

\end{proof}

\end{document}